\documentclass[11pt,reqno]{amsart}

\usepackage{geometry}
\usepackage{amsmath,amssymb,amsfonts,mathrsfs}

\usepackage{graphicx}

\usepackage[non-sorted-cites,initials,alphabetic,nobysame]{amsrefs}

\usepackage{url}
\usepackage{comment}

\usepackage{amsthm}
\usepackage{cleveref}
\usepackage{aliascnt}
\theoremstyle{definition}

\newtheorem{theorem}{Theorem}[section]
\crefname{theorem}{Theorem}{Theorems}
\Crefname{theorem}{Theorem}{Theorems}

\newaliascnt{definition}{theorem}

\aliascntresetthe{definition}
\crefname{definition}{Definition}{Definitions}
\Crefname{definition}{Definition}{Definitions}

\newaliascnt{proposition}{theorem}
\newtheorem{proposition}[proposition]{Proposition}
\aliascntresetthe{proposition}
\crefname{proposition}{Proposition}{Propositions}
\Crefname{proposition}{Proposition}{Propositions}

\newaliascnt{lemma}{theorem}
\newtheorem{lemma}[lemma]{Lemma}
\aliascntresetthe{lemma}
\crefname{lemma}{Lemma}{Lemmas}
\Crefname{lemma}{Lemma}{Lemmas}

\newaliascnt{corollary}{theorem}

\aliascntresetthe{corollary}
\crefname{corollary}{Corollary}{Corollaries}
\Crefname{corollary}{Corollary}{Corollaries}

\newaliascnt{remark}{theorem}
\newtheorem{remark}[remark]{Remark}
\aliascntresetthe{remark}
\crefname{remark}{Remark}{Remarks}
\Crefname{remark}{Remark}{Remarks}

\newaliascnt{example}{theorem}

\aliascntresetthe{example}
\crefname{example}{Example}{Examples}
\Crefname{example}{Example}{Examples}

\newaliascnt{assumption}{theorem}

\aliascntresetthe{assumption}
\crefname{assumption}{Assumption}{Assumptions}
\Crefname{assumption}{Assumption}{Assumptions}

\numberwithin{equation}{section}
\crefname{equation}{equation}{equations}
\Crefname{equation}{Equation}{Equations}
\crefname{figure}{figure}{figures}
\Crefname{figure}{Figure}{Figures}
\crefname{table}{table}{tables}
\Crefname{table}{Table}{Tables}

\crefname{section}{Section}{Sections}
\Crefname{section}{Section}{Sections}
\crefname{subsection}{Section}{Sections}
\Crefname{subsection}{Section}{Sections}
\crefname{appendix}{Appendix}{Appendices}
\Crefname{appendix}{Appendix}{Appendices}

\DeclareMathOperator*{\argmin}{arg\,min} 

\DeclareMathOperator{\Exp}{Exp}

\newcommand{\dual}{*}

\newcommand{\energy}{\mathrm{E}}
\newcommand{\length}{\mathrm{L}}

\newcommand{\curves}{\Xi}

\begin{document}

\title[Minimization of Curve Length through Energy Minimization]{
    Minimization of Curve Length through Energy Minimization \\
    using Finite Differences and Numerical Integration \\
    in Euclidean Space
}
\author{
    Akira Kitaoka
}
\address{NEC Corporation, 1753 Shimonumabe, Nakahara-ku, Kawasaki, Kanagawa, Japan }
\email{akira-kitaoka@nec.com}

\thanks{
    We would like to thank Rin Takano, Taichi Hirano, Yuzuru Okajima, and Yoichi Sasaki for their careful reading of this paper.
}


\begin{abstract}
    We consider the approximation of minimal geodesics between two closed sets in $\mathbb{R}^D$ endowed with a smooth Riemannian metric. The continuous problem is formulated as the minimization of the energy functional over piecewise smooth curves joining the two sets. We study discrete approximations obtained by finite differences together with numerical integration, and reconstruct continuous curves from discrete minimizers by linear interpolation.
    Our main result is a direct convergence analysis of the trapezoidal-rule discretization. We prove that the energy of the linearly interpolated discrete minimizer converges to the minimum energy, and that the squared length of the reconstructed curve converges to the squared minimal length, both with rate $O(N^{-1/2})$ as the number of subintervals $N$ tends to infinity. We also obtain the corresponding reconstruction estimates for the left-endpoint rule. In addition, we give an explicit example showing that a direct discretization of the length functional does not, in general, converge to the minimal length.

\end{abstract}


\maketitle

\section{Introduction}

The problem of determining minimal geodesics, namely shortest paths with respect to a given metric, arises in various areas of mathematics and its applications. Examples include path detection and shortest-route computation \cites{caselles1997geodesic,cohen1997global,deschamps2001fast}, segmentation \cites{bai2007geodesic,price2010geodesic}, animation \cite{heeren2012time}, meshing \cite{liu2017constructing}, shape exploration \cite{rabin2010geodesic}, optimal trajectory planning in robotics \cites{zhang2007robot,zhang2010optimal}, satellite orbits \cites{ghaffari1970integration}, clothing design \cites{sanchez-reyes2008constrainted}, and geodesic deep learning \cites{masci2015geodesic,he2019geonet}. For Riemannian manifolds, existing approximation methods for minimal geodesics and their lengths include mesh-based reductions to graph shortest path problems \cites{cohen1997global,deschamps2001fast}, shooting methods based on the geodesic equation (cf. \Cref{eq:def_geodesic}) \cite{patrikalakis2003computational}, boundary value methods for the geodesic equation \cites{patrikalakis2003computational,kasap2005numerical,seif2018numerical}, and discrete variational formulations of energy minimization using finite differences and numerical integration \cites{heeren2012time,rumpf2013discrete,rumpf2015variational,meng20211variational,liu2023efficiently}.

A standard variational approach is to approximate the energy functional by finite differences and numerical integration, and then to solve the resulting finite-dimensional optimization problem. Compared with mesh-based methods, this avoids discretizing the ambient space itself, which may be costly in high dimensions. Compared with formulations based directly on the geodesic equation, it uses only the Riemannian metric and does not require Christoffel symbols explicitly. In the fixed-endpoint setting, \cite{rumpf2015variational} established $O(N^{-1})$ control of the discrete minimum value for consistent discretizations of the energy on manifolds embedded in Euclidean space. The present paper is concerned with different quantities: the energy and squared length of the continuous curves reconstructed from discrete minimizers by linear interpolation, and with the set-to-set setting in which the endpoints are allowed to vary in prescribed closed sets.

In this paper, we study the approximation of energy minimization between two closed sets in $\mathbb{R}^D$ endowed with a smooth Riemannian metric. Our main focus is on the trapezoidal-rule discretization and on the continuous curves reconstructed from discrete minimizers by linear interpolation. We prove that the energy of the reconstructed curve converges to the minimum energy, and that the squared length of the reconstructed curve converges to the squared minimal length, both with rate $O(N^{-1/2})$ as the number of subintervals $N$ tends to infinity. We also treat the left-endpoint rule in the same set-to-set framework, mainly for comparison and for the analysis of reconstructed curves. For the left-endpoint rule, the rate obtained here is not claimed to be optimal. 
In particular, the $O(N^{-1})$ estimate from \cite{rumpf2015variational} concerns the discrete minimum value, whereas our main reconstruction estimates concern the continuous energy and squared length of the linearly interpolated discrete minimizers. Thus the two types of estimates are complementary rather than contradictory. We also give an explicit rigorous example showing that a direct discretization of the length functional does not, in general, recover the minimal length, even as $N \to \infty$.

Before stating the main theorems, we introduce some notation. Let $D \in \mathbb{Z}_{\geq 1}$ be an integer, $p \in [1 , \infty ]$, and $\| \bullet \|_p$ denote the $\ell^p$ norm in the real coordinate space $\mathbb{R}^D$. Consider non-empty sets $X^{(0)} , X^{(1)} \subset \mathbb{R}^D$. In the main theorems below, we further assume that they are closed and that at least one of them is bounded. When $M$ is a manifold, we denote by $\curves_{\mathrm{ps}}$ the set of continuous piecewise $C^1$ curves $\gamma \colon [0,1] \to M$. For a set of curves $\curves \subset \curves_{\mathrm{ps}} $, we define
\[ 
    \curves (X^{(0)} , X^{(1)} ) := \{ \gamma \in \curves \, | \, \gamma (0) \in X^{(0)} , \, \gamma (1) \in X^{(1)} \}
    .
\]
When it is clear from the context, we identify $\{ x \}$ with the point $x \in \mathbb{R}^D$. Let $N \in \mathbb{Z}_{\geq 1}$ be a positive integer, $n = 0,\ldots, N$, and set constants $t_n = n/N$, $h = 1/N$. The set of point sequences in $\mathbb{R}^D$ is
\[
    \curves_{\mathrm{p}}^N 
    := \{ \gamma^{\mathrm{p}} \colon \{ t_n | n = 0, \ldots,N \} \to \mathbb{R}^D \}
\]
and we define
\[
    \curves_{\mathrm{p}}^N (X^{(0)} , X^{(1)}) := 
    \{ \gamma^{\mathrm{p}} \in \curves_{\mathrm{p}}^N \, | \, \gamma^{\mathrm{p}} (0) \in X^{(0)}, \,
    \gamma^{\mathrm{p}} (1) \in X^{(1)}
    \}
\]
as well.
Finite differences are defined by
\begin{equation*}
    \beta^{\mathrm{p}} (t_n) := \frac{1}{h} (\gamma^{\mathrm{p}} (t_{n+1})  - \gamma^{\mathrm{p}} (t_{n}) ), \, \text{ for } n =0 , \ldots , N-1
\end{equation*}
and this corresponds to the derivative $\dot{\gamma}$ for a curve $\gamma \in \curves_{\mathrm{ps}}$. For a point sequence\\ $\gamma^{\mathrm{p}} \in \curves_{\mathrm{p}}^N (X^{(0)} , X^{(1)})$, the linear interpolation $\gamma^{\mathrm{pl}}$ is defined as follows: for any $n = 0, \ldots , N-1$, when $0 \leq s \leq h$, writing $t = t_n + s$, we have
\[
    \gamma^{\mathrm{pl}} (t) :=
        (1 - N s ) \gamma^{\mathrm{p}} (t_n) + Ns \gamma^{\mathrm{p}} (t_{n+1}) 
        .
\]
The energy $\energy^{g} ( \gamma)$ and length $\length^g ( \gamma)$ of a curve $\gamma \in \curves_{\mathrm{ps}}$ with respect to the Riemannian metric $g$ are defined as follows:
\begin{equation*}
    \energy^g ( \gamma) := \int_{0}^{1} g_{\gamma (t)}  \left( \dot{\gamma} (t) , \dot{\gamma} (t) \right) dt ,
    \quad
    \length^g ( \gamma) := \int_{0}^{1} g_{\gamma (t)} \left( \dot{\gamma} (t) , \dot{\gamma} (t) \right)^{1/2} dt .
\end{equation*}
The distance between two points $x^{(0)},x^{(1)} \in \mathbb{R}^D$ is defined by
\[
    d^g (x^{(0)},x^{(1)}) := \inf_{\gamma \in \curves_{\mathrm{ps}} (x^{(0)},x^{(1)})} \length^g ( \gamma),
\]
and for sets $X^{(0)}, X^{(1)} \subset \mathbb{R}^D$, it is defined by
\[
    d^g (X^{(0)},X^{(1)})
    :=
    \inf_{\substack{
        x^{(0)} \in X^{(0)} ,\\
        x^{(1)} \in X^{(1)}
    } }
    d^g (x^{(0)},x^{(1)})
\]
The energy $\energy^{g}$, approximated by finite differences and numerical integration (trapezoidal rule), is defined for any sequence of points $\gamma^{\mathrm{p}} \colon \{ t_n \} \to \mathbb{R}^D$ as follows:
\begin{align*}
    \energy^{g}_{\mathrm{tra},N} (\gamma^{\mathrm{p}})
    :=
    \frac{1}{N} \sum_{n=0}^{N-1} 
            g_{\frac{\gamma^{\mathrm{p}} (t_n ) + \gamma^{\mathrm{p}} (t_{n+1} )}{2}}  ( \beta^{\mathrm{p}} (t_n), \beta^{\mathrm{p}} (t_{n}) )
    .
\end{align*}
The energy $\energy^{g}$, approximated by finite differences and numerical integration (left-endpoint rule), is defined for any $\gamma^{\mathrm{p}} \in \curves_{\mathrm{p}}^N$ as follows:
\begin{equation*}
    \energy^{g}_{\mathrm{l},N} (\gamma^{\mathrm{p}})
    :=
    \frac{1}{N} \sum_{n=0}^{N-1} 
            g_{\gamma^{\mathrm{p}} (t_n)} ( \beta^{\mathrm{p}} (t_n), \beta^{\mathrm{p}} (t_{n}) )
    .
\end{equation*}
The point sequences $\gamma_{\mathrm{tra},N}^{*\mathrm{p}}$ and $\gamma_{\mathrm{l},N}^{*\mathrm{p}}$ are chosen as minimizers:
\begin{align}
    \gamma_{\mathrm{tra},N}^{*\mathrm{p}} 
    & \in 
    \argmin_{
        \gamma^{\mathrm{p}} 
        \in 
        \curves_{\mathrm{p}}^N (X^{(0)}, X^{(1)}) 
    } 
        \energy_{\mathrm{tra},N}^{g} ( \gamma^{\mathrm{p}} )
    ,
    \label{eq:approx_numerical_integral_trapezoidal}
    \\ 
    \gamma_{\mathrm{l},N}^{*\mathrm{p}} 
    & \in 
    \argmin_{
        \gamma^{\mathrm{p}} 
        \in 
        \curves_{\mathrm{p}}^N (X^{(0)}, X^{(1)}) 
    } 
        \energy_{\mathrm{l},N}^{g} ( \gamma^{\mathrm{p}} )
    .
    \label{eq:approx_numerical_integral_left}
\end{align}
The linear interpolations of the point sequences $\gamma_{\mathrm{tra},N}^{*\mathrm{p}}$ and $\gamma_{\mathrm{l},N}^{*\mathrm{p}}$ are denoted by $\gamma_{\mathrm{tra},N}^{*\mathrm{pl}}$ and $\gamma_{\mathrm{l},N}^{*\mathrm{pl}}$, respectively.
Their existence is proved later in \Cref{prop:existence_discrete_minimizers}.

The error when applying the trapezoidal rule to the energy is as follows.
\begin{theorem}
    \label{theo:approx_numerical_integral_trapezoidal}
    Let $g \in C^{\infty}$ be a Riemannian metric on $\mathbb{R}^D$ such that there exist constants $c_1,c_2>0$ satisfying
    \[
    c_1 \|u\|_2^2 \le g_x(u,u) \le c_2 \|u\|_2^2
    \quad\text{for all } x,u\in\mathbb{R}^D,
    \]
    and such that, writing $g_x(u,u)=u^\top H(x)u$, one has
    \[
    \|H(x)-H(y)\|_{\mathcal B}\le L_H\|x-y\|_2
    \quad\text{for all } x,y\in\mathbb{R}^D.
    \]
    We assume that $X^{(0)}, X^{(1)}$ are closed sets in $\mathbb{R}^D$, and suppose that at least one of them is bounded.
    Then, there exists a constant $C$ which depends only on $c_1, c_2, L_H, d^g (X^{(0)}, X^{(1)}), X^{(0)}, X^{(1)}$ such that
    \begin{align}
        \begin{split}
        & \min_{\gamma \in \curves_{\mathrm{ps}} (X^{(0)}, X^{(1)}) } \energy^{g}( \gamma )
        -
        \frac{C }{N^{1/2}}
        \leq 
        \energy^{g} ( \gamma_{\mathrm{tra},N}^{*\mathrm{pl}}  )
        -
        \frac{C }{N^{1/2}}
        \\
        & \leq 
        \energy^{g}_{\mathrm{tra},N} ( \gamma_{\mathrm{tra},N}^{*\mathrm{p}}  )
        \leq  
        \min_{\gamma \in \curves_{\mathrm{ps}} (X^{(0)}, X^{(1)}) } \energy^{g}( \gamma )
        +
        \frac{C }{N}.
        \end{split}
        \notag
    \end{align}
\end{theorem}
\Cref{theo:approx_numerical_integral_trapezoidal} shows that the continuous curve obtained by linearly interpolating the discrete minimizer for the trapezoidal rule has energy within $O(N^{-1/2})$ of the minimum energy. It also provides the corresponding upper estimate for the discrete minimum value.

The key to the proof of \Cref{theo:approx_numerical_integral_trapezoidal} involves the following steps. First, for any given curve, the error between the energy approximation using the trapezoidal rule and the true energy is evaluated by applying Morrey-type inequalities to the $L^2$ norms of the first- and second-order derivatives of the curve. Next, the $L^2$ norms of the first- and second-order derivatives of a geodesic are evaluated. Furthermore, when the curve is close to the minimum energy, it belongs to a neighborhood of the geodesic, allowing evaluation of the $L^2$ norms of the first- and second-order derivatives of all curves belonging to this neighborhood. Finally, for the neighborhood of the geodesic, the error between the energy approximated using the trapezoidal rule and the true energy is evaluated, thereby proving \Cref{theo:approx_numerical_integral_trapezoidal}.

The squared length of the reconstructed curve $\gamma_{\mathrm{tra},N}^{*\mathrm{pl}}$ converges to the squared minimal length at the rate $O(N^{-1/2})$.
\begin{theorem}
    \label{theo:approx_numerical_integral_trapezoidal_length}
    We assume that the setting of \Cref{theo:approx_numerical_integral_trapezoidal}.
    Then, there exists a constant $C$ which depends only on $c_1, c_2, L_H, d^g (X^{(0)}, X^{(1)}), X^{(0)}, X^{(1)}$ such that
    \begin{align} 
        \length^{g} ( \gamma_{\mathrm{tra},N}^{*\mathrm{pl}}  )^2
        -
        \min_{\gamma \in \curves_{\mathrm{ps}} (X^{(0)}, X^{(1)}) } \length^{g}( \gamma )^2
        \leq
        \frac{C }{N^{1/2}}.
        \notag
    \end{align}
\end{theorem}

We next consider the left-endpoint rule in the same set-to-set framework, mainly for comparison with the trapezoidal rule and for the analysis of reconstructed curves.
The error when applying the left-endpoint rule to the energy is as follows.
\begin{theorem}
    \label{theo:approx_numerical_integral_left}
    We assume that the setting of \Cref{theo:approx_numerical_integral_trapezoidal}.
    Then, there exists a constant $C$ which depends only on $c_1, c_2, L_H, d^g (X^{(0)}, X^{(1)}), X^{(0)}, X^{(1)}$ such that
    \begin{align}
        \begin{split}
        & \min_{\gamma \in \curves_{\mathrm{ps}} (X^{(0)}, X^{(1)}) } \energy^{g}( \gamma )
        -
        \frac{C }{N^{1/2}}
        \leq 
        \energy^{g} (\gamma_{\mathrm{l},N}^{*\mathrm{pl}} )
        -
        \frac{C }{N^{1/2}}
        \\
        & \leq 
        \energy^{g}_{\mathrm{l},N} ( \gamma_{\mathrm{l},N}^{*\mathrm{p}} )
        \leq  
        \min_{\gamma \in \curves_{\mathrm{ps}} (X^{(0)}, X^{(1)}) } \energy^{g}( \gamma )
        +
        \frac{C }{N^{1/2}}.
        \end{split}
        \notag
    \end{align}
\end{theorem}
\Cref{theo:approx_numerical_integral_left} provides the analogous reconstruction estimate for the left-endpoint rule in the same set-to-set setting.

The key to the proof of \Cref{theo:approx_numerical_integral_left} involves the following steps. First, for any given curve, the error between the energy approximations using the trapezoidal rule and the left-endpoint rule is evaluated using the $L^2$ norm of the finite differences. Using this estimate, the proof of \Cref{theo:approx_numerical_integral_left} proceeds in essentially the same way as that of \Cref{theo:approx_numerical_integral_trapezoidal}.

The squared length of the reconstructed curve $\gamma_{\mathrm{l},N}^{*\mathrm{pl}}$ converges to the squared minimal length at the rate $O(N^{-1/2})$.
\begin{theorem}
    \label{theo:approx_numerical_integral_left_length}
    We assume that the setting of \Cref{theo:approx_numerical_integral_trapezoidal}.
    Then, there exists a constant $C$ which depends only on $c_1, c_2, L_H, d^g (X^{(0)}, X^{(1)}), X^{(0)}, X^{(1)}$ such that
    \begin{align} 
        \length^{g} ( \gamma_{\mathrm{l},N}^{*\mathrm{pl}}  )^2
        -
        \min_{\gamma \in \curves_{\mathrm{ps}} (X^{(0)}, X^{(1)}) } \length^{g}( \gamma )^2
        \leq
        \frac{C }{N^{1/2}}.
        \notag
    \end{align}
\end{theorem}
\cite[Lemma 4.7]{rumpf2015variational} shows that $\mathcal{W}_{\mathrm{l}}[\mathbf{x},\tilde{\mathbf{x}}]$, corresponding to the left-endpoint rule, fits into the framework of \cite{rumpf2015variational}. The trapezoidal integrand
\[
    \mathcal{W}_{\mathrm{tra}}[\mathbf{x},\tilde{\mathbf{x}}]
    :=
    g_{(\mathbf{x}+\tilde{\mathbf{x}})/2}(\tilde{\mathbf{x}}-\mathbf{x},\tilde{\mathbf{x}}-\mathbf{x})
\]
also appears to satisfy the same third-order consistency requirement used there. However, the reconstruction estimates for the linearly interpolated discrete minimizers obtained in the present paper, as well as the set-to-set formulation treated here, are not provided in that form by \cite{rumpf2015variational}. For this reason, we give a direct analysis of the trapezoidal discretization in $\mathbb{R}^D$, starting from the metric $g$ itself.

From the computational point of view, both discretizations lead to finite-dimensional optimization problems; see, for example, \cites{xie2017differential,jallet2022constrained,rosenbrock1972differential} for related optimization methods. We do not pursue algorithmic aspects in the present paper.

The structure of this paper is organized as follows. In \Cref{sec:related_work}, we describe the related research. In \Cref{sec:dg}, we recall differential geometry. In \Cref{sec:minimal_geodesic}, we show that the minimal curve is essentially the minimal geodesic. \Cref{sec:length_minimization_and_numerical_integral} provides an example where the minimization problem approximating length via numerical integration fails to approach the minimal length, irrespective of the number of points $N$. In \Cref{sec:length_numerical_integral_min}, we show the equivalence between the energy minimization problem and finding the minimal geodesic for complete Riemannian manifolds. In \Cref{sec:proof_of_approx_numerical_integral} we prove  \Cref{theo:approx_numerical_integral_trapezoidal}. In \Cref{sec:proof_of_approx_trapezoidal_left}, we prove  \Cref{theo:approx_numerical_integral_left}. Finally, in \Cref{sec:proof_of_length_approx_trapezoidal_left} we prove \Cref{theo:approx_numerical_integral_trapezoidal_length,theo:approx_numerical_integral_left_length}.

\section{Related work}
\label{sec:related_work}

\subsection*{Minimization of Lagrangian functionals and numerical integration}
Examples of approximating functionals by finite differences and numerical integration can be found in works such as \cites{futurologist2020,hong2024convergence}. Theorem 3.2 in \cite{hong2024convergence} shows that, for the Freidlin--Wentzell action functional, the error between the minimum of the functional and the minimum of its discrete approximation is $O(N^{-1/2})$.
\Cref{theo:approx_numerical_integral_trapezoidal,theo:approx_numerical_integral_left} treat the case in which the functional is the energy $\energy^{g}$.

\subsection*{Length minimization problem and numerical integration}

To approximate the length of a minimal geodesic, one may consider solving a nonlinear optimization problem obtained by approximating the length functional via numerical integration.
However, \cite[Remark 4.12]{rumpf2015variational} sketches a counterexample showing that the length obtained through this method does not asymptotically approach the minimal path length regardless of the number of points $N$ used, but does not provide a rigorous non-convergence proof. We explicitly construct a non-convergent example (\Cref{sec:length_numerical_integral_min}). 
This indicates that solving an optimization problem obtained by numerically approximating the length functional does not, in general, recover the desired minimal path.

\subsection*{Formulation of the energy minimization problem}

In addition to the trapezoidal rule and the left-endpoint rule \cite{meng20211variational}, \cite{liu2017constructing} provides another way to approximately formulate the energy minimization problem by numerical integration. \cite{liu2017constructing} uses the fact that, along a geodesic, the Riemannian norm of the velocity is constant, and formulates the problem as follows:
\begin{align}
    \begin{split}
        \text{minimize }&
        \energy^{g}_{\mathrm{l},N} (\gamma^{\mathrm{p}})
        \\
        \text{subject to } &
            g_{\gamma^{\mathrm{p}} (t_n)} ( \beta^{\mathrm{p}} (t_n), \beta^{\mathrm{p}} (t_{n}) )
            =
            g_{\gamma^{\mathrm{p}} (t_{n+1})} ( \beta^{\mathrm{p}} (t_{n+1}), \beta^{\mathrm{p}} (t_{n+1}) )
            \text{ for } n = 0, \ldots , N-2.
    \end{split}
    \label{eq:liu2017}
\end{align}
Removing the constraints in \Cref{eq:liu2017} yields the formulation proposed in \cites{heeren2012time,rumpf2015variational}. \cite{rumpf2015variational} showed that, when $M$ is a manifold embedded in Euclidean space equipped with a smooth Riemannian metric $g$ and $X^{(0)}, X^{(1)}$ are singleton sets, one has
\begin{align}
    \left|
        \energy^{g}_{\mathrm{l},N} ( \gamma_{\mathrm{l},N}^{*\mathrm{p}} )
        -\min_{\gamma \in \curves_{\mathrm{ps}} (X^{(0)}, X^{(1)}) } \energy^{g}( \gamma )
    \right|
    \leq
    \frac{C }{N}
    \notag
\end{align}
by \cite[Theorem 4.5, Lemma 4.7]{rumpf2015variational}.
For the convenience of the reader, we summarize the relation between \cite{rumpf2015variational} and the present paper in \Cref{tab:comparison_rumpf_wirth}.

\begin{table}[t]
    \centering
    \small
    \begin{tabular}{|p{3.2cm}|p{4.6cm}|p{4.6cm}|}
        \hline
        Aspect & \cite{rumpf2015variational} & Present paper \\
        \hline
        Endpoint condition & fixed endpoints & set-to-set endpoints \\
        \hline
        Discrete minimum value & treated; $O(N^{-1})$ for consistent discretizations & treated for comparison \\
        \hline
        Energy of linearly interpolated discrete minimizers & not treated explicitly & treated; $O(N^{-1/2})$ reconstruction estimate \\
        \hline
        Squared length of reconstructed curves & not treated explicitly & treated; $O(N^{-1/2})$ estimate \\
        \hline
        Direct discretization of length & counterexample sketched in Remark 4.12 & explicit rigorous non-convergence example \\
        \hline
    \end{tabular}
    \caption{Comparison with \cite{rumpf2015variational}.}
    \label{tab:comparison_rumpf_wirth}
\end{table}

\section{Differential geometry}
\label{sec:dg}

For any $p \in [1, \infty)$ and the interval $[s_1, s_2] \subset \mathbb{R}$, we define the norm $\|\cdot\|_{L^p([s_1, s_2])}$ for any (Lebesgue measurable) mapping $\gamma \colon [s_1, s_2] \to \mathbb{R}^D$ as follows:
\[
    \| \gamma \|_{L^p([s_1, s_2])}^p := \int_{s_1}^{s_2} \left\| \gamma (t) \right\|_2^p \, dt.
\]
When the context is clear, we may simplify $L^p([s_1, s_2])$ as $L^p$. For a matrix $A \in \mathbb{R}^{D \times D}$, we define the operator norm as $\| A \|_{\mathcal{B}} := \sup_{x \in \mathbb{R}^D \setminus \{0\}} \| A x \|_2 / \|x \|_2$.

Let $M$ be a manifold. Denote the tangent bundle of $M$ by $TM$ and the cotangent bundle by $T^{\dual}M$. A function $g_{\bullet} (\bullet , \bullet ) \in C^2 ( M, T^{\dual}M \otimes T^{\dual}M)$ is a Riemannian metric on $M$ if, for any $x \in M$, $g_x$ defines a real inner product on $T_x M$. The pair $(M,g)$ is called a Riemannian manifold. When clear from the context, we may denote $g(u,v)$ as $g_x(u,v)$. A specific example of a Riemannian metric is given by taking $M = \mathbb{R}^D$ and specifying a positive definite matrix $H(x)$ for each $x \in M$ such that, for any $u, v \in \mathbb{R}^D$, we have $g_x(u,v) = u^{\top} H(x) v$. In this case, the pair $(\mathbb{R}^D , g)$ forms a Riemannian manifold. Additionally, the Euclidean metric is denoted by $g_{\mathrm{Euc},x} (u,u) := \|u \|_2^2$.

In a local neighborhood $U$ with local coordinates $x=(x_1, \ldots, x_D)$, 
\[
    g_{ij} (x) = g_x \left( 
        \frac{\partial}{\partial x_i}, 
        \frac{\partial}{\partial x_j} 
    \right)
\]
are of class $C^2$. At each point $x \in U$, the matrix $g = (g_{ij})$ is positive definite, which allows the inverse matrix $g^{-1}=(g^{ij})$ to be defined. The Christoffel symbols of the second kind for the Levi-Civita connection, $\Gamma_{ij}^k \in C^1$, are given by
\begin{equation}
    \Gamma_{i j}^k (x)
    = \frac{1}{2} \sum_{\ell = 1}^{D} g^{k \ell} \left(
        \frac{\partial g_{\ell j}}{\partial x_i} 
        + 
        \frac{\partial g_{i \ell}}{\partial x_j}
        -
        \frac{\partial g_{i j}}{\partial x_\ell}
    \right).
    \label{eq:Christoffel_symbol}
\end{equation}
A path $\gamma \in \curves_{\mathrm{ps}} \cap C^{2}$ is called a geodesic if, in each local coordinate neighborhood $U$ through which the curve $\gamma$ passes, it satisfies the geodesic equation:
\begin{equation}
    \ddot{\gamma}_k = -\sum_{i,j} \Gamma_{i j}^k (\gamma) \dot{\gamma}_i \dot{\gamma}_j.
    \label{eq:def_geodesic}
\end{equation}
The set of all geodesics is denoted by $\curves_{\mathrm{geo}}$.

For a connected Riemannian manifold $(M, g)$, the distance $d^g$ between any two points $x^{(0)}, x^{(1)} \in M$ is defined by
\[
    d^g (x^{(0)}, x^{(1)}) := \inf_{\gamma \in \curves_{\mathrm{ps}} (x^{(0)}, x^{(1)} )} \length^g ( \gamma).
\]
A path $\gamma \in \curves_{\mathrm{ps}}(x^{(0)}, x^{(1)})$ is said to be minimal if it satisfies
\[
    d^g (x^{(0)}, x^{(1)}) = \length^g ( \gamma).
\]

Henceforth, unless otherwise specified, we assume that the metric $g$ is smooth (i.e., $g \in C^{\infty}$).
\begin{proposition}[Cf. \cite{milnor1963morse}*{\S 10}]
    \label{prop:metric_geodesic_constant}
    Let $(M, g)$ be a Riemannian manifold. Then, for any geodesic $\gamma \in \curves_{\mathrm{geo}}$, 
    \[
        g_{\gamma (t)} (\dot{\gamma}(t) , \dot{\gamma}(t))
    \]
    is constant for $t \in [0,1]$. 
    In particular, we have
    \[
        g_{\gamma (t)} (\dot{\gamma}(t) , \dot{\gamma}(t)) = \length^g ( \gamma)^2.
    \]
\end{proposition}

A connected Riemannian manifold $(M, g)$ is called complete if the metric space $(M, d^g)$ is complete.

\begin{proposition}[\cite{hopf1931ueber}]
    \label{prop:Hopf-Rinow-theorem}
    Let $(M, g)$ be a connected and complete Riemannian manifold. Then, there exists a minimal geodesic connecting any two points in $M$.
\end{proposition}

A sufficient condition for completeness is given in the following proposition.
\begin{proposition}
    \label{prop:bounded_g_completeness}
    Let $g$ be a Riemannian metric on $\mathbb{R}^D$ such that there exist constants $c_1, c_2 > 0$ for which $c_1 \|u \|_2^2 \leq g_x (u,u) \leq c_2 \| u \|_2^2$ holds for all $x, u \in \mathbb{R}^D$. Then, the manifold $(\mathbb{R}^D, g)$ is complete.
\end{proposition}

\begin{proof}
    For any curve $\gamma \in \curves_{\mathrm{ps}}$, we have
    \[
        c_1^{1/2} \length^{g_{\mathrm{Euc}}} ( \gamma)
        \leq 
        \length^{g} ( \gamma)
        \leq 
        c_2^{1/2} \length^{g_{\mathrm{Euc}}} ( \gamma).
    \]
    Applying $\inf_{\gamma \in \curves_{\mathrm{ps}} (x^{(0)}, x^{(1)} )}$ to this inequality, we obtain
    \[
        c_1 \| x^{(0)} - x^{(1)} \|^2
        \leq 
        d^g ( x^{(0)}, x^{(1)} )^2
        \leq
        c_2 \| x^{(0)} - x^{(1)} \|^2.
    \]
    Since the Euclidean distance is complete on $\mathbb{R}^D$, $(\mathbb{R}^D, g)$ is also complete.
\end{proof}

\section{Minimal curves and minimal geodesics}
\label{sec:minimal_geodesic}

In this section, we show that a path of minimal length is essentially a minimal geodesic. For any $x \in M$, let $V_x$ be an open neighborhood of $0$ in the tangent space $T_x M$, and define $\Exp_x \colon V_x \to M$ as the exponential map on the Riemannian manifold $(M, g)$ with respect to the Levi-Civita connection $\nabla$. 

For any $x \in M$, if there exists an open neighborhood $V_x$ of $0$ in $T_x M$ such that $\Exp_x \colon V_x \to \Exp_x (V_x)$ is a diffeomorphism, then $V_x$ is called a normal coordinate neighborhood of $x$, and $\Exp_x (V_x)$ is referred to as a normal neighborhood of $x$.

Geodesics are known to be locally length-minimizing.
\begin{proposition}[Cf. \cite{konno2013differentiable}*{Theorem 4.2.1}]
    \label{prop:geodesic_min_loc_path_min}
    Let $(M, g)$ be a Riemannian manifold. Let $U$ be a normal neighborhood of $x^{(0)} \in M$, with some $\varepsilon > 0$ such that $B = \{x \in M \mid d^g (x, x^{(0)}) < \varepsilon \} \subset U$. Let $x^{(1)} \in B$. Then, a curve $\gamma^* \in \curves_{\mathrm{ps}}(x^{(0)}, x^{(1)})$ satisfies
    \[
        \length^g (\gamma^*) \leq 
        \length^g (\gamma)
    \]
    for any curve $\gamma \in \curves_{\mathrm{ps}}(x^{(0)}, x^{(1)})$ if and only if there exists a minimal geodesic $\gamma^{**} \in \curves_{\mathrm{geo}} (x^{(0)}, x^{(1)})$, and a monotone increasing continuous piecewise $C^1$ function $\phi \colon [0, 1] \to [0, 1]$ with $\phi(0) = 0$, $\phi(1) = 1$, such that for any $t \in [0, 1]$, we have
    \[
        \gamma^{*} (t) = \gamma^{**} \circ \phi (t).
    \]
\end{proposition}

We now establish the following characterization.
\begin{proposition}
    \label{prop:geodesic_min_path_min}
    Let $(M, g)$ be a Riemannian manifold. Let $x^{(0)}, x^{(1)} \in M$ be arbitrary points. Then, a curve $\gamma^* \in \curves_{\mathrm{ps}}(x^{(0)}, x^{(1)})$ satisfies
    \[
        \length^g (\gamma^*) \leq 
        \length^g (\gamma)
    \]
    for any curve $\gamma \in \curves_{\mathrm{ps}}(x^{(0)}, x^{(1)})$ if and only if there exists a minimal geodesic $\gamma^{**} \in \curves_{\mathrm{geo}} (x^{(0)}, x^{(1)})$, and a monotone increasing continuous piecewise $C^1$ function $\phi \colon [0, 1] \to [0, 1]$ with $\phi(0) = 0$, $\phi(1) = 1$, such that for any $t \in [0, 1]$, we have
    \[
        \gamma^{*} (t) = \gamma^{**} \circ \phi (t).
    \]
\end{proposition}

Before proving this proposition, we recall how the exponential map provides a canonical local coordinate system.

\begin{proposition}[Cf. \cite{konno2013differentiable}*{Theorem 4.2.1}]
    \label{prop:uniformal_normal_radius}
    Let $(M, g)$ be a Riemannian manifold. Let $K \subset M$ be a compact set. Then, there exists some $\varepsilon > 0$ such that for any $x \in K$, the exponential map $\Exp_x \colon \{ v \in T_x M \mid g_x (v, v)^{1/2} < \varepsilon \} \to \{ y \in M \mid d^g (x, y ) < \varepsilon \}$ is a diffeomorphism.
\end{proposition}

\begin{proof}[Proof of \Cref{prop:geodesic_min_path_min}]

    Step 1: 
    For any $[ t^{(0)} , t^{(1)} ] \subset [0,1]$, let $\gamma \in \curves_{\mathrm{ps}} (\gamma^* (t^{(0)}), \gamma^* (t^{(1)}))$. Then, we have
    \begin{equation}
        \int_{t^{(0)}}^{t^{(1)}} 
            g_{\gamma^* (t)} \left( \dot{\gamma}^* (t) , \dot{\gamma}^* (t) \right)^{1/2} dt
        \leq 
        \length^g (\gamma).
        \label{eq:local_short_path}
    \end{equation}
    Suppose, for contradiction, that there exists a curve $\gamma^{(1)} \colon [ t^{(0)} , t^{(1)} ] \to M$ such that
    \begin{equation}
        \int_{t^{(0)}}^{t^{(1)}} 
            g_{\gamma^{(1)} (t)} \left( \dot{\gamma}^{(1)} (t) , \dot{\gamma}^{(1)} (t) \right)^{1/2} dt
        <
        \int_{t^{(0)}}^{t^{(1)}} 
            g_{\gamma^* (t)} \left( \dot{\gamma}^* (t) , \dot{\gamma}^* (t) \right)^{1/2} dt.
    \end{equation}
    Define the curve $\tilde{\gamma}^{(1)} \in \curves_{\mathrm{ps}}(x^{(0)}, x^{(1)})$ by
    \begin{equation}
    \tilde{\gamma}^{(1)} (t)
    :=
        \begin{cases}
            \gamma^{(1)} (t), & t \in [ t^{(0)} , t^{(1)} ], \\
            \gamma^* (t), & \text{otherwise}.
        \end{cases}
        \notag  
    \end{equation}
    Then, $\length^g (\tilde{\gamma}^{(1)} ) < \length^g (\gamma^* )$, resulting in a contradiction.

    Step 2:
    Let $t^{(0)} \in [0,1]$. Consider a normal neighborhood $U = \left\{ x \in M \middle| d^g (\gamma^* (t^{(0)}), x) < \varepsilon \right\}$ of the point $\gamma^* (t^{(0)})$. Since the curve $\gamma^*$ is continuous, there exists an interval $[t^{(1)}, t^{(2)}] \ni t^{(0)}$ such that
    \[
        \gamma^* ([ t^{(1)} ,t^{(2)} ]) \subset U,
        \quad
        d^g (\gamma^* (t^{(0)}), \gamma^* (t^{(1)})) = \frac{\varepsilon}{2},
        \quad
        d^g (\gamma^* (t^{(0)}), \gamma^* (t^{(2)})) = \frac{\varepsilon}{2}.
    \]
    By \Cref{prop:geodesic_min_loc_path_min}, there exist a geodesic $\gamma^{**}_2 \in \curves_{\mathrm{geo}}$ with a monotone increasing continuous piecewise $C^1$ function $\phi_2 \colon [t^{(1)}, t^{(0)}] \to [0,1]$ satisfying $\phi_2 (t^{(1)} ) = 0$, $\phi_2 (t^{(0)}) = 1$, and a geodesic $\gamma^{**}_3 \in \curves_{\mathrm{geo}}$ with a monotone increasing continuous piecewise $C^1$ function $\phi_3 \colon [t^{(0)}, t^{(2)}] \to [0,1]$ satisfying $\phi_3 (t^{(0)} ) = 0$, $\phi_3 (t^{(2)}) = 1$, such that
    \begin{equation}
        \gamma^* (t) = \gamma^{**}_2 \circ \phi_2 (t), \text{ for } t \in [t^{(1)}, t^{(0)}],
        \quad
        \gamma^* (t) = \gamma^{**}_3 \circ \phi_3 (t), \text{ for } t \in [t^{(0)}, t^{(2)}].
        \label{eq:local_geodesic_min_loc_path_min}
    \end{equation}

    Step 3:
    By applying Proposition \ref{prop:uniformal_normal_radius} to the set $\gamma^* ([0,1])$, there exists some $\varepsilon > 0$ such that for any $x \in \gamma^* ([0,1])$, the exponential map $\Exp_x \colon \{ v \in T_x M \mid g_x (v, v)^{1/2} < \varepsilon \} \to \{ y \in M \mid d^g (x, y) < \varepsilon \}$ is a diffeomorphism. Therefore, $\{ y \in M \mid d^g (x, y) < \varepsilon \}$ is a normal coordinate neighborhood of $x \in \gamma^* ([0,1])$.

    Step 4: We show the necessity of this condition.
    Let $\{ \gamma^* (t^{(n)}) \}_{n=1}^N$ be a sequence of points such that for any $n = 1, \ldots , N-1$
    \begin{equation}
        d^g (\gamma^* (t^{(n)}), \gamma^* (t^{(n+1)})) < \frac{\varepsilon}{4},
        \label{eq:near_points}
    \end{equation}
    and for $i \neq j$, 
    \begin{equation}
        \gamma^* (t^{(i)}) \neq \gamma^* (t^{(j)}).
        \label{eq:different_point}
    \end{equation}
    For any $n = 1, \ldots , N-2$, we have $d^g (\gamma^* (t^{(n)}), \gamma^* (t^{(n+2)})) < \varepsilon /2$.
    According to \Cref{eq:local_geodesic_min_loc_path_min}, there exists a geodesic $\gamma^{*n} \in \curves_{\mathrm{geo}}$, and a monotone increasing continuous piecewise $C^1$ function $\phi_n \colon [t^{(n)}, t^{(n+2)}] \to [0,1]$ with $\phi_n (t^{(n)}) = 0$, $\phi_n (t^{(n+2)}) = 1$ such that
    \[
        \gamma^* (t) = \gamma^{*n} \circ \phi_n (t), \quad t \in [t^{(n)}, t^{(n+2)}].
    \]
    Since geodesics locally satisfy \Cref{eq:def_geodesic}, the uniqueness of solutions to the geodesic equation implies that there exist a geodesic $\gamma^{**}$ and a monotone increasing continuous piecewise $C^1$ function $\phi \colon [0,1] \to [0,1]$ with $\phi(0) = 0$, $\phi(1) = 1$ satisfying
    \[
        \gamma^* = \gamma^{**} \circ \phi.
    \]

    Step 5: We show the sufficiency of the condition.

    Let $\{ \gamma^* (t^{(n)}) \}_{n=1}^N$ be a sequence of points satisfying \Cref{eq:near_points,eq:different_point} for all $n = 1, \ldots , N-1$. By \Cref{eq:local_geodesic_min_loc_path_min} and \Cref{prop:geodesic_min_loc_path_min}, we have
    \begin{equation}
        \int_{t^{(n)}}^{t^{(n+1)}} 
            g_{\gamma^{*} (t)} \left( \dot{\gamma}^{*} (t) , \dot{\gamma}^{*} (t) \right)^{1/2} dt
        =
        \int_{\gamma^{**-1}(\gamma^* (t^{(n)}))}^{\gamma^{**-1}(\gamma^* (t^{(n+1)}))} 
            g_{\gamma^{**} (t)} \left( \dot{\gamma}^{**} (t) , \dot{\gamma}^{**} (t) \right)^{1/2} dt
            .
    \end{equation}
    By summing over $n$, we obtain
    \[
        \length^g (\gamma^{*})
        =
        \length^g (\gamma^{**}).
    \]
    Since $\gamma^{**}$ is a minimal geodesic, this proves sufficiency.
\end{proof}

\begin{lemma}
    \label{lem:minimal_constant_speed_geodesic}
    Let $(M,g)$ be a Riemannian manifold, and let
    $\gamma \in \curves_{\mathrm{ps}}(x^{(0)},x^{(1)})$ be a minimal path.
    Assume further that
    \[
        g_{\gamma(t)}(\dot{\gamma}(t),\dot{\gamma}(t))
    \]
    is constant on $[0,1]$.
    Then $\gamma$ is a geodesic.
\end{lemma}

\begin{proof}
    By \Cref{prop:geodesic_min_path_min}, there exist a minimal geodesic
    $\eta \in \curves_{\mathrm{geo}}(x^{(0)},x^{(1)})$ and a monotone
    increasing continuous piecewise $C^1$ map
    $\phi \colon [0,1] \to [0,1]$
    satisfying
    \[
        \phi(0)=0,\quad \phi(1)=1,
    \]
    such that
    \[
        \gamma = \eta \circ \phi.
    \]

    Since $\eta$ is a geodesic, \Cref{prop:metric_geodesic_constant} yields that
    \[
        g_{\eta(s)}(\dot{\eta}(s),\dot{\eta}(s))
    \]
    is constant on $[0,1]$. Hence
    \[
        g_{\gamma(t)}(\dot{\gamma}(t),\dot{\gamma}(t))
        =
        g_{\eta(\phi(t))}(\dot{\eta}(\phi(t)),\dot{\eta}(\phi(t)))\,|\phi'(t)|^2
    \]
    and the assumption that it is constant imply that
    $|\phi'(t)|$ is constant almost everywhere on $[0,1]$.

    Furthermore, since $\phi$ is monotone increasing and satisfies
    \[
        \phi(0)=0,\quad \phi(1)=1,
    \]
    it follows that
    \[
        \phi(t)=t.
    \]
    Therefore
    \[
        \gamma=\eta,
    \]
    and $\gamma$ is a geodesic.
\end{proof}

\section{Minimal length problem and numerical integration}
\label{sec:length_minimization_and_numerical_integral}

In this section, we present an example showing that a minimization problem based on a numerical approximation of the length functional does not approximate the true length minimization problem, even if the number of points $N$ is increased.
For this section only, consider the function $f(x_1, x_2) := -\cos x_1 + 2$, with the metric $g^f := f^2 g_{\mathrm{Euc}}$, where $g_{\mathrm{Euc}}$ is the Euclidean metric. The points are given by $x^{(0)} = (0,0)$ and $x^{(1)} = (4\pi,0)$. For this section only, we set $X^{(0)} = \{ x^{(0)} \} $ and $X^{(1)} = \{ x^{(1)} \} $.

\subsection{Length of minimal path}

Geometrically, the conformal factor
\[
    f(x_1,x_2)=2-\cos x_1
\]
creates ``valleys'' at $x_1=0,2\pi,4\pi,\ldots$, where the metric is smallest. The true minimal curve from $(0,0)$ to $(4\pi,0)$ must nevertheless traverse the whole interval in the $x_1$-direction, and its weighted length is $8\pi$. By contrast, a direct discretization of the length functional samples the metric only at finitely many points, so the discrete optimizer can jump between such low-cost locations and artificially reduce the discrete length to $4\pi$, independently of $N$.

To compute the minimal path length, observe that if a path $\gamma$ is minimal from $x^{(0)}$ to $x^{(1)}$ in the metric $g^f$, then $\gamma_2(t) = 0$. This follows since
\begin{align}
    \length^{g^f} (\gamma) 
    &= \int_0^1 \left( -\cos(\gamma_1(t)) + 2 \right) \sqrt{(\dot{\gamma}_1(t))^2 + (\dot{\gamma}_2(t))^2} \, dt
    \notag \\
    &\geq \int_0^1 \left( -\cos(\gamma_1(t)) + 2 \right) |\dot{\gamma}_1(t)| \, dt
    = \length^{g^f} \left((\gamma_1, 0)\right). 
    \notag
\end{align}
Additionally, if $\gamma$ is a minimal path from $x^{(0)}$ to $x^{(1)}$ in the metric $g^f$, then $\gamma_1(t)$ must be monotone increasing in $t$. This can be proved by contradiction. Suppose there exists an interval $[t_1, t_2]$ over which $\gamma_1(t)$ is strictly decreasing. Then either $\gamma_1(t_1) > 0$ or $\gamma_1(t_2) < 4\pi$ must be satisfied.
Assume first that $\gamma_1(t_1) > 0$. By the Intermediate Value Theorem, there exists some $t_3 \in [0, t_1]$ such that $\gamma_1(t_3) = \gamma_1(t_2)$. Define a new path

\[
    \check{\gamma}(t) := 
    \begin{cases}
        \gamma(t_2), & t_3 \leq t \leq t_2, \\
        \gamma(t), & \text{otherwise}.
    \end{cases}
\]
This yields $\length^{g^f} (\gamma) > \length^{g^f} (\check{\gamma})$, which contradicts the assumption that $\gamma$ is a minimal path. Similarly, if $\gamma_1(t_2) < 4\pi$, a contradiction arises through a similar argument. Thus, $\gamma_1(t)$ must be monotone increasing in $t$.

Furthermore, we have $\gamma ([0,1]) = [0,4 \pi] \times \{ 0 \}$. This is because, by definition, $\gamma(0) = (0,0)$ and $\gamma(1) = (4\pi,0)$, and since $[0,1]$ is path-connected, $\gamma ([0,1])$ is also path-connected. Since $\gamma_2(t) = 0$, it follows that $\gamma ([0,1]) \supset [0,4 \pi] \times \{ 0 \}$.

On the other hand, suppose there exists a point $(x_1, 0) \in \gamma ([0,1])$ with $x_1 \not\in [0,4\pi]$. If $x_1 < 0$, there would exist some $t_3 \in [0,1]$ such that $\gamma_1 (t_3) = x_1 < 0$. By the Intermediate Value Theorem, there must be some $t_2 \in (t_3, 1]$ such that $\gamma_1 (t_2) = 0$. Define a new path by

\[
    \check{\gamma} (t) := 
    \begin{cases}
        \gamma (0) , &  0 \leq t \leq t_2, \\
        \gamma (t) , & \text{otherwise}.
    \end{cases}
\]

This results in $\length^{g^f} (\gamma) > \length^{g^f} (\check{\gamma})$, which contradicts the assumption that $\gamma$ is a minimal path. A similar contradiction arises if $x_1 > 4\pi$. Therefore, $\gamma ([0,1]) = [0,4 \pi] \times \{ 0 \}$.

To sum up, if a curve $\gamma \in \curves_{\mathrm{ps}} (x^{(0)}, x^{(1)})$ is a minimal path, then $\gamma_2(t)=0$ and $\gamma_1 \colon [0,1] \to [0,4 \pi]$ is monotone increasing with $\gamma_1(0) = 0$ and $\gamma_1(1) = 4\pi$. Additionally, since $(\mathbb{R}^2, g^f)$ is complete, there exists a minimal geodesic $\gamma^*$ from $x^{(0)}$ to $x^{(1)}$ with respect to $g^f$. This geodesic is unique. Solving \Cref{eq:def_geodesic}, we find that $\gamma_1^* \colon [0,1] \to [0,4\pi]$ is strictly increasing and is a diffeomorphism.

To show that $\hat{\gamma}(t) = t(4\pi,0)$ is a minimal curve, define $\phi := \gamma_1^{*-1} \circ \hat{\gamma}_1$. Here, $\phi$ is monotone increasing and satisfies $\phi(0)=0$ and $\phi(1)=1$. Therefore, it follows that $\hat{\gamma} = \gamma^* \circ \phi$. According to \Cref{prop:geodesic_min_path_min}, $\hat{\gamma}$ is indeed minimal.

The length of $\hat{\gamma}$ with respect to $g^f$ is computed as follows:
\begin{align}
    \length^{g^f} (\hat{\gamma})
    &= \int_0^1 \left(-\cos(8 \pi t) + 2\right) \cdot 4 \pi \, dt
    = 8 \pi. 
    \notag
\end{align}
Thus, we have
\begin{equation}
    \min_{\gamma \in \curves_{\mathrm{ps}} (x^{(0)}, x^{(1)}) }
    \length^{g^f} (\gamma )
    = 8 \pi.
    \label{eq:length-min}
\end{equation}

\subsection{Left-endpoint rule}

We approximate the length of a curve using the left-endpoint rule as follows:
\[
    \length^{g}_{\mathrm{l},N} (\gamma^{\mathrm{p}})
    :=
    \frac{1}{N} \sum_{n=0}^{N-1} 
            g_{\gamma^{\mathrm{p}} (t_n)} (\beta^{\mathrm{p}} (t_n), \beta^{\mathrm{p}} (t_n))^{1/2}.
\]
We show that the minimization problem for $\length^{g}_{\mathrm{l},N}$ does not coincide with the true length minimization problem, regardless of the choice of $N$.

The approximation of length using the left-endpoint rule can be evaluated as follows:
\begin{align}
    \length^{g}_{\mathrm{l},N} (\gamma^{\mathrm{p}})
    & =
    \frac{1}{N} \sum_{n=0}^{N-1} 
            g_{\gamma^{\mathrm{p}} (t_n)} ( \beta^{\mathrm{p}} (t_n), \beta^{\mathrm{p}} (t_{n}) )^{1/2}
    \notag \\
    & =
    \frac{1}{N} \sum_{n=0}^{N-1} 
            f (\gamma^{\mathrm{p}} (t_n)) \| \beta^{\mathrm{p}} (t_n) \|_2
    \notag \\
    & \geq 
    \sum_{n=0}^{N-1} 
            \| 
                \gamma^{\mathrm{p}} (t_{n+1}) 
                -
                \gamma^{\mathrm{p}} (t_n) 
            \|_2
    & \because f(x_1,x_2) \geq 1
    \notag \\
    & \geq 
    \| 
                \gamma^{\mathrm{p}} (t_{N}) 
                -
                \gamma^{\mathrm{p}} (t_0) 
    \|_2
    = 4 \pi
    \notag 
    & \because \text{triangle inequality}.
\end{align}
We have
\begin{equation}
    \min_{
        \gamma^{\mathrm{p}} 
        \in 
        \curves_{\mathrm{p}}^N (X^{(0)}, X^{(1)}) 
    } 
        \length_{\mathrm{l},N}^{g} ( \gamma^{\mathrm{p}} )
    \geq 4 \pi.
    \label{eq:length-left-point-min-bound}
\end{equation}

On the other hand, consider a sequence of points $\{\gamma^{*\mathrm{p}} (t_n)\}_{n=0}^N$ such that
\[
    \gamma_2^{*\mathrm{p}}(t_n)=0,
    \quad
    \gamma_1^{*\mathrm{p}} (t_n) \in \{0, 2\pi, 4\pi\},
    \quad
    \gamma_1^{*\mathrm{p}} (t_n) \leq \gamma_1^{*\mathrm{p}} (t_{n+1}),
    \gamma^{*\mathrm{p}} \in \curves_{\mathrm{p}}^N (X^{(0)}, X^{(1)})
\]
In this case, since $f (\gamma^{*\mathrm{p}} (t_n)) = 1$, we obtain
\begin{equation}
    \length^{g}_{\mathrm{l},N} (\gamma^{*\mathrm{p}})
    = 4 \pi.
    \label{eq:length-left-point-example-min}
\end{equation}
From \Cref{eq:length-left-point-min-bound,eq:length-left-point-example-min}, we have
\begin{equation}
    \min_{
        \gamma^{\mathrm{p}} 
        \in 
        \curves_{\mathrm{p}}^N (X^{(0)}, X^{(1)}) 
    } 
        \length_{\mathrm{l},N}^{g} ( \gamma^{\mathrm{p}} )
    = 4 \pi.
    \label{eq:length-left-point-min}
\end{equation}
From \Cref{eq:length-min,eq:length-left-point-min}, we conclude that the approximation of length using the left-endpoint rule does not approximate the minimal length, regardless of how large $N$ is chosen.

\subsection{Trapezoidal Rule}

We approximate the length of a curve using the trapezoidal rule as follows:
\[
    \length^{g}_{\mathrm{tra},N} (\gamma^{\mathrm{p}})
    :=
    \frac{1}{N} \sum_{n=0}^{N-1} 
            g_{\frac{\gamma^{\mathrm{p}} (t_n) + \gamma^{\mathrm{p}} (t_{n+1})}{2}} (\beta^{\mathrm{p}} (t_n), \beta^{\mathrm{p}} (t_n))^{1/2}.
\]
We show that the minimization problem for $\length^{g}_{\mathrm{tra},N}$ does not coincide with the true length minimization problem, regardless of how large $N$ is chosen.

The approximation of length using the trapezoidal rule can be evaluated as follows:
\begin{align}
    \length^{g}_{\mathrm{tra},N} (\gamma^{\mathrm{p}})
    & = 
    \frac{1}{N} \sum_{n=0}^{N-1} 
            g_{\frac{\gamma^{\mathrm{p}} (t_n ) + \gamma^{\mathrm{p}} (t_{n+1} )}{2}} ( \beta^{\mathrm{p}} (t_n), \beta^{\mathrm{p}} (t_{n}) )^{1/2}
        \notag \\
    & = 
    \frac{1}{N} \sum_{n=0}^{N-1} 
            f \left(
                \frac{\gamma^{\mathrm{p}} (t_n ) + \gamma^{\mathrm{p}} (t_{n+1} )}{2}
            \right) 
            \| \beta^{\mathrm{p}} (t_n) \|_2
    \notag \\
        & \geq 
    \sum_{n=0}^{N-1} 
            \| 
                \gamma^{\mathrm{p}} (t_{n+1}) 
                -
                \gamma^{\mathrm{p}} (t_n) 
            \|_2
    & \because f(x_1,x_2) \geq 1
    \notag \\
    & \geq 
    \| 
                \gamma^{\mathrm{p}} (t_{N}) 
                -
                \gamma^{\mathrm{p}} (t_0) 
    \|_2
    = 4 \pi
    & \because \text{triangle inequality}.
    \notag 
\end{align}
Thus, we have
\begin{equation}
    \min_{
        \gamma^{\mathrm{p}} 
        \in 
        \curves_{\mathrm{p}}^N (X^{(0)}, X^{(1)}) 
    } 
        \length_{\mathrm{tra},N}^{g} ( \gamma^{\mathrm{p}} )
    \geq 4 \pi.
    \label{eq:length-tra-point-min-bound}
\end{equation}
On the other hand, consider a sequence of points $\{\gamma^{*\mathrm{p}} (t_n)\}_{n=0}^N$ such that
\[
    \gamma_2^{*\mathrm{p}}(t_n)=0,
    \quad
    \gamma_1^{*\mathrm{p}} (t_n) \in \{0, 4\pi\},
    \quad
    \gamma_1^{*\mathrm{p}} (t_n) \leq \gamma_1^{*\mathrm{p}} (t_{n+1}),
    \quad
    \gamma^{*\mathrm{p}} \in \curves_{\mathrm{p}}^N (X^{(0)}, X^{(1)}).
\]
Since
\[
    f \left(
                \frac{\gamma^{*\mathrm{p}} (t_n ) + \gamma^{*\mathrm{p}} (t_{n+1} )}{2}
            \right)  
            =1 ,
\]
we obtain
\begin{equation}
    \length^{g}_{\mathrm{tra},N} (\gamma^{*\mathrm{p}})
    = 4 \pi .
    \label{eq:length-tra-point-example-min}
\end{equation}

From \Cref{eq:length-tra-point-min-bound,eq:length-tra-point-example-min}, we have
\begin{equation}
    \min_{
        \gamma^{\mathrm{p}} 
        \in 
        \curves_{\mathrm{p}}^N (X^{(0)}, X^{(1)}) 
    } 
        \length_{\mathrm{tra},N}^{g} ( \gamma^{\mathrm{p}} )
    = 4 \pi
    .
    \label{eq:length-tra-point-min}
\end{equation}
From \Cref{eq:length-min,eq:length-tra-point-min}, we conclude that the approximation of length using the trapezoidal rule does not approximate the minimal length, regardless of how large $N$ is chosen.

\section{Characterization of minimal geodesics via energy minimization}
\label{sec:length_numerical_integral_min}

In this section, we show that for a complete Riemannian manifold, the problem of minimizing energy between two sets is equivalent to finding the minimal geodesic between them.

\subsection{Case of two points}

We recall the equivalence between the energy minimization problem and the problem of finding the minimal geodesic between two points.
\begin{proposition}[\cite{milnor1963morse}*{Corollary 10.7, Lemma 12.1}]
    \label{prop:geodesic_len_and_energy}
    Let $(M, g)$ be a Riemannian manifold. Assume that there exists a minimal geodesic in $(M, g)$ connecting the points $x^{(0)}, x^{(1)} \in M$. Then, we have
    \begin{equation*}
        \min_{\gamma \in \curves_{\mathrm{ps}} (x^{(0)}, x^{(1)})} \energy^g (\gamma)
        = 
        \min_{\gamma \in \curves_{\mathrm{ps}} (x^{(0)}, x^{(1)})} \length^g (\gamma)^2,
        \quad
        \argmin_{\gamma \in \curves_{\mathrm{ps}} (x^{(0)}, x^{(1)})} \energy^g (\gamma)
        = 
        \argmin_{\gamma \in \curves_{\mathrm{geo}} (x^{(0)}, x^{(1)})} \length^g (\gamma).
    \end{equation*}
    In particular, 
    \[
        \argmin_{\gamma \in \curves_{\mathrm{ps}} (x^{(0)}, x^{(1)})} \energy^g (\gamma)
    \]
    consists of geodesics.
\end{proposition}

\subsection{Case of two sets}

We show the equivalence between the energy minimization problem and the problem of finding the minimal geodesic between two sets.

\begin{proposition}
    \label{prop:bounded_2sets_min}
    Let $(M,g)$ be a complete Riemannian manifold, and let $X^{(0)}, X^{(1)} \subset M$ be non-empty closed sets, at least one of which is bounded. Then, the set
    \[
         \mathcal{Y} := \{ (x^{(0)}, x^{(1)}) \in X^{(0)} \times X^{(1)} \, | \, d^g (x^{(0)}, x^{(1)}) = d^g (X^{(0)}, X^{(1)}) \}
    \]
    is non-empty, bounded, and closed.
\end{proposition}

\begin{proof}
    By symmetry, assume $X^{(0)}$ is bounded. Since $X^{(1)}$ is non-empty, we can choose an element $x^{(1)} \in X^{(1)}$. Define the constant $R$ and the set $X^{(1)\prime}$ as follows:
    \begin{align}
        R & := \max_{x^{(0)} \in X^{(0)}} d^g (x^{(0)}, x^{(1)}), \notag \\
        X^{(1)\prime} & := \{ x^{(1)} \in X^{(1)} \mid d^g (X^{(0)}, x^{(1)}) \leq R \}.
        \notag
    \end{align}
    Since $X^{(0)}$ is bounded, the set $X^{(1)\prime}$ is bounded and closed. Moreover, we have
    \begin{equation}
        d^g (X^{(0)}, X^{(1)}) = d^g (X^{(0)}, X^{(1)\prime}).
    \end{equation}
    The function $d^g \colon M \times M \to \mathbb{R}$ is continuous. By a standard consequence of \Cref{prop:Hopf-Rinow-theorem}, on a complete Riemannian manifold every closed bounded set is compact. Hence both $X^{(0)}$ and $X^{(1)\prime}$ are compact, so $X^{(0)} \times X^{(1)\prime}$ is compact as well. Since $d^g \colon M \times M \to \mathbb{R}$ is continuous, it attains its minimum on $X^{(0)} \times X^{(1)\prime}$. Therefore $\mathcal{Y}$ is non-empty.

    The set $\mathcal{Y}$ is the preimage of the value $d^g (X^{(0)}, X^{(1)})$ under the continuous function $d^g$. Therefore $\mathcal{Y}$ is closed. Furthermore, since $\mathcal{Y} \subset X^{(0)} \times X^{(1)\prime}$ and $X^{(0)} \times X^{(1)\prime}$ is bounded, $\mathcal{Y}$ is also bounded.
\end{proof}

\begin{proposition}
    \label{prop:geodesic_len_and_energy_2}
    Let $(M,g)$ be a complete Riemannian manifold, and let $X^{(0)}, X^{(1)} \subset M$ be non-empty closed sets, at least one of which is bounded. Then, we have
    \begin{align*}
        \min_{\gamma \in \curves_{\mathrm{ps}} (X^{(0)} , X^{(1)} )} \energy^g ( \gamma)
        & = 
        \min_{\gamma \in \curves_{\mathrm{ps}} (X^{(0)} , X^{(1)} )} \length^g ( \gamma)^2,
        \\
        \argmin_{\gamma \in \curves_{\mathrm{ps}} (X^{(0)} , X^{(1)} )} \energy^g ( \gamma)
        & = 
        \argmin_{\gamma \in \curves_{\mathrm{geo}} (X^{(0)} , X^{(1)} )} \length^g ( \gamma) \not = \emptyset.
    \end{align*}
    In particular, 
    \[
        \argmin_{\gamma \in \curves_{\mathrm{ps}} (X^{(0)} , X^{(1)} )} \energy^g ( \gamma)
    \]
    consists of geodesics.
\end{proposition}

\begin{proof}
    By symmetry, we assume that $X^{(0)}$ is bounded. We first show that
    \begin{equation}
        \argmin_{\gamma \in \curves_{\mathrm{ps}} (X^{(0)} , X^{(1)} )} \length^g ( \gamma) \not = \emptyset,
        \quad
        d^g (X^{(0)} , X^{(1)}) =
        \min_{\gamma \in \curves_{\mathrm{ps}} (X^{(0)} , X^{(1)} )} \length^g ( \gamma)
        =
        \min_{\gamma \in \curves_{\mathrm{geo}} (X^{(0)} , X^{(1)} )} \length^g ( \gamma)
        \label{eq:argmin_length_path}
    \end{equation}
    holds.

    By \Cref{prop:bounded_2sets_min}, there exist $x^{(0)} \in X^{(0)}$ and $x^{(1)} \in X^{(1)}$ such that
    \[
        d^g(x^{(0)},x^{(1)}) = d^g(X^{(0)},X^{(1)}).
    \]
    Moreover, by \Cref{prop:Hopf-Rinow-theorem}, there exists a minimal geodesic connecting $x^{(0)}$ and $x^{(1)}$. Therefore,
    \[
        d^g (X^{(0)} , X^{(1)})
        =
        \min_{\gamma \in \curves_{\mathrm{ps}} (X^{(0)} , X^{(1)} )} \length^g ( \gamma)
        =
        \min_{\gamma \in \curves_{\mathrm{geo}} (X^{(0)} , X^{(1)} )} \length^g ( \gamma),
    \]
    which yields \Cref{eq:argmin_length_path}.

    Next, we show
    \begin{equation}
        \min_{\gamma \in \curves_{\mathrm{ps}} (X^{(0)} , X^{(1)} )} \energy^g ( \gamma)
        =
        \min_{\gamma \in \curves_{\mathrm{ps}} (X^{(0)} , X^{(1)} )} \length^g ( \gamma)^2
        \label{eq:argmin_length_path_01}
    \end{equation}
    and
    \begin{equation}
        \argmin_{\gamma \in \curves_{\mathrm{ps}} (X^{(0)} , X^{(1)} )} \energy^g ( \gamma)
        =
        \argmin_{\gamma \in \curves_{\mathrm{geo}} (X^{(0)} , X^{(1)} )} \length^g ( \gamma).
        \label{eq:argmin_length_path_02}
    \end{equation}

    By the Cauchy--Schwarz inequality, for any $\gamma \in \curves_{\mathrm{ps}}$,
    \begin{equation}
        \length^g(\gamma)^2 \leq \energy^g(\gamma),
        \label{eq:CS-length-energy-fixed}
    \end{equation}
    and equality is equivalent to
    \[
        g_{\gamma(t)}(\dot{\gamma}(t),\dot{\gamma}(t))
    \]
    being constant almost everywhere on $[0,1]$.

    By \Cref{eq:argmin_length_path}, for any $\gamma \in \curves_{\mathrm{ps}}(X^{(0)},X^{(1)})$,
    \[
        \energy^g(\gamma)^{1/2}
        \geq
        \length^g(\gamma)
        \geq
        d^g(X^{(0)},X^{(1)}).
    \]
    Hence, equality holds if and only if $\gamma$ is a minimal curve between $X^{(0)}$ and $X^{(1)}$ and
    \[
        g_{\gamma(t)}(\dot{\gamma}(t),\dot{\gamma}(t))
    \]
    is constant.

    By \Cref{lem:minimal_constant_speed_geodesic}, this is in turn equivalent to $\gamma$ being a minimal geodesic connecting some\\
    $(x^{(0)},x^{(1)}) \in X^{(0)} \times X^{(1)}$ with
    \[
        d^g(x^{(0)},x^{(1)}) = d^g(X^{(0)},X^{(1)}).
    \]

    Therefore, \Cref{eq:argmin_length_path_01,eq:argmin_length_path_02} follow.
\end{proof}

\section{Approximation of energy minimization problem via trapezoidal rule}
\label{sec:proof_of_approx_numerical_integral}

In this section, we prove \Cref{theo:approx_numerical_integral_trapezoidal}. In \Cref{sec:Morrey_type_ineq}, we prove a Morrey-type inequality. In \Cref{sec:error_of_energy_and_trapezoidal_rule}, we estimate the error between the energy and its trapezoidal approximation when an $L^2$ upper bound exists for first and second derivatives in the space of curves. 
In \Cref{sec:trapolized_and_interpolation}, we estimate the error between the discrete energy defined by the trapezoidal rule and the energy of the curve obtained by linearly interpolating the corresponding point sequence.
In \Cref{sec:estimatation_Christoffel_symbol}, we evaluate the first derivative of the metric and the Christoffel symbols. In \Cref{sec:apriori_estimate_for_minimal_geodesic}, we estimate the $L^2$ norms of the first and second derivatives of minimal geodesics. Finally, in \Cref{sec:proof_of_theo:approx_numerical_integral_trapezoidal}, we prove \Cref{theo:approx_numerical_integral_trapezoidal}.

We first record the existence of discrete minimizers for the trapezoidal-rule
and left-endpoint-rule energies.

\begin{proposition}
    \label{prop:existence_discrete_minimizers}
    Let $g \in C^0$ be a Riemannian metric on $\mathbb{R}^D$ such that
    \[
        c_1 \|u\|_2^2 \le g_x(u,u)
        \quad\text{for all } x,u \in \mathbb{R}^D
    \]
    with some constant $c_1>0$.
    Let $X^{(0)}, X^{(1)} \subset \mathbb{R}^D$ be non-empty closed sets, and suppose that at least one of them is bounded.
    Then, for every $N \in \mathbb{Z}_{\ge 1}$, the sets
    \[
        \argmin_{\gamma^{\mathrm{p}} \in \curves_{\mathrm{p}}^N(X^{(0)},X^{(1)})}
        \energy_{\mathrm{tra},N}^g(\gamma^{\mathrm{p}})
        \quad\text{and}\quad
        \argmin_{\gamma^{\mathrm{p}} \in \curves_{\mathrm{p}}^N(X^{(0)},X^{(1)})}
        \energy_{\mathrm{l},N}^g(\gamma^{\mathrm{p}})
    \]
    are non-empty.
\end{proposition}

\begin{proof}
    By symmetry, it suffices to treat the case where $X^{(0)}$ is bounded.

    Identify a point sequence $\gamma^{\mathrm{p}} \in \curves_{\mathrm{p}}^N$ with the tuple
    \[
        (x_0,\ldots,x_N)
        :=
        (\gamma^{\mathrm{p}}(t_0),\ldots,\gamma^{\mathrm{p}}(t_N))
        \in (\mathbb{R}^D)^{N+1}.
    \]
    Then the admissible set
    \[
        \mathcal{X}_N
        :=
        \{(x_0,\ldots,x_N)\in(\mathbb{R}^D)^{N+1}
        \mid x_0\in X^{(0)},\ x_N\in X^{(1)}\}
    \]
    is closed.

    Since $g$ is continuous, both functionals
    \[
        (x_0,\ldots,x_N)\mapsto \energy_{\mathrm{tra},N}^g(\gamma^{\mathrm{p}})
        \quad\text{and}\quad
        (x_0,\ldots,x_N)\mapsto \energy_{\mathrm{l},N}^g(\gamma^{\mathrm{p}})
    \]
    are continuous on $(\mathbb{R}^D)^{N+1}$.

    Let $R>0$, and let $\gamma^{\mathrm{p}} \in \curves_{\mathrm{p}}^N(X^{(0)},X^{(1)})$ satisfy
    \[
        \energy_{\mathrm{tra},N}^g(\gamma^{\mathrm{p}}) \le R.
    \]
    Writing $x_n := \gamma^{\mathrm{p}}(t_n)$, we have
    \[
        \energy_{\mathrm{tra},N}^g(\gamma^{\mathrm{p}})
        =
        \frac1N \sum_{n=0}^{N-1}
        g_{\frac{x_n+x_{n+1}}2}\bigl(\beta^{\mathrm{p}}(t_n),\beta^{\mathrm{p}}(t_n)\bigr)
        \ge
        \frac{c_1}{N}\sum_{n=0}^{N-1}\|\beta^{\mathrm{p}}(t_n)\|_2^2.
    \]
    Since $\beta^{\mathrm{p}}(t_n)=N(x_{n+1}-x_n)$, this yields
    \[
        c_1 N \sum_{n=0}^{N-1}\|x_{n+1}-x_n\|_2^2
        \le R.
    \]
    Hence, by the Cauchy--Schwarz inequality,
    \[
        \sum_{n=0}^{N-1}\|x_{n+1}-x_n\|_2
        \le
        N^{1/2}
        \left(\sum_{n=0}^{N-1}\|x_{n+1}-x_n\|_2^2\right)^{1/2}
        \le
        \left(\frac{R}{c_1}\right)^{1/2}.
    \]
    Therefore, for every $n=0,\ldots,N$,
    \[
        \|x_n\|_2
        \le
        \|x_0\|_2 + \sum_{k=0}^{N-1}\|x_{k+1}-x_k\|_2
        \le
        \sup_{x\in X^{(0)}}\|x\|_2 + \left(\frac{R}{c_1}\right)^{1/2}.
    \]
    Thus
    \[
        \{\gamma^{\mathrm{p}} \in \curves_{\mathrm{p}}^N(X^{(0)},X^{(1)})
        \mid \energy_{\mathrm{tra},N}^g(\gamma^{\mathrm{p}})\le R\}
    \]
    is bounded. Since it is also closed, it is compact.

    By continuity of $\energy_{\mathrm{tra},N}^g$, the minimum of $\energy_{\mathrm{tra},N}^g$ over
    $\curves_{\mathrm{p}}^N(X^{(0)},X^{(1)})$ is attained.

    Since
    \[
        \energy_{\mathrm{l},N}^g(\gamma^{\mathrm{p}})
        =
        \frac1N \sum_{n=0}^{N-1}
        g_{x_n}\bigl(\beta^{\mathrm{p}}(t_n),\beta^{\mathrm{p}}(t_n)\bigr)
        \ge
        \frac{c_1}{N}\sum_{n=0}^{N-1}\|\beta^{\mathrm{p}}(t_n)\|_2^2,
    \]
    in the same way as $\energy_{\mathrm{tra},N}^g$,
    we get the proof for $\energy_{\mathrm{l},N}^g$.
\end{proof}

\subsection{Morrey-type Inequality}
\label{sec:Morrey_type_ineq}

In this subsection, we recall Morrey's inequality and derive a related estimate from its proof.
\begin{proposition}[Morrey Inequality (cf. \cite{evans2010partial}*{\S 5.6.2, Theorem 4})]
    \label{prop:morrey_ineq}
    Let $\gamma \in \curves_{\mathrm{ps}}$ satisfy $\| \gamma \|_{L^2}, \| \dot{\gamma} \|_{L^2} < \infty$. Then, for any $0 \leq s_1 \leq s \leq s_2 \leq 1$, we have
    \begin{align}
        \| \gamma (s_2) - \gamma (s_1 ) \|_2 
        & \leq (s_2 - s_1)^{1/2} \| \dot{\gamma} \|_{L^2([s_1,s_2])} ,
        \label{eq:Morrey_ineq_1}\\
        \| \gamma (s) \|_2 
        & \leq 
        \| \gamma \|_{L^2([s_1,s_2])}
        + 2 \| \dot{\gamma} \|_{L^2([s_1,s_2])} .
        \label{eq:Morrey_ineq_2}
    \end{align}
\end{proposition}

Moreover, a slight modification of the proof of Morrey's inequality yields the following estimate.
\begin{proposition}
    \label{prop:morrey_ineq_cor}
    Let $\gamma \in \curves_{\mathrm{ps}}$ satisfy $\| \dot{\gamma} \|_{L^2}, \| \ddot{\gamma} \|_{L^2} < \infty$. Then, for any $0 \leq s_1 \leq s \leq s_2 \leq 1$, we have
    \begin{equation}
        \| \gamma (s_2) - \gamma (s_1 ) \|_2 
        \leq (s_2 - s_1) \left( \| \dot{\gamma} \|_{L^2([s_1,s_2])}
        + 2 \| \ddot{\gamma} \|_{L^2([s_1,s_2])} \right).
        \notag
    \end{equation}
\end{proposition}

\begin{proof}
    By the fundamental theorem of calculus, we have
    \begin{equation}
        \| \gamma (s_2) - \gamma (s_1 ) \|_2
        \leq \left\| \int_{s_1}^{s_2} \dot{\gamma}(s) ds \right\|_2
        \leq \int_{s_1}^{s_2} \| \dot{\gamma}(s) \|_2 ds.
        \notag
    \end{equation}    
    From \Cref{eq:Morrey_ineq_2}, we have
    \begin{align}
        \int_{s_1}^{s_2} \| \dot{\gamma}(s) \|_2 ds
        & \leq 
        \int_{s_1}^{s_2} \left( 
        \| \dot{\gamma} \|_{L^2([s_1,s_2])}
        + 2 \| \ddot{\gamma} \|_{L^2([s_1,s_2])} 
        \right)
        ds \notag \\ 
        & = (s_2 - s_1) \left( \| \dot{\gamma} \|_{L^2([s_1,s_2])}
        + 2 \| \ddot{\gamma} \|_{L^2([s_1,s_2])} \right).
        \notag 
    \end{align}
\end{proof}

For a curve $\gamma \in \curves_{\mathrm{ps}}$, we also define the discrete difference at each $n = 0, \ldots, N-1$ by
\[
    \beta_{\gamma}(t_n)
    :=
    \frac{\gamma(t_{n+1}) - \gamma(t_n)}{h}.
\]
When the context is clear, we write $\beta(t_n)$ instead of $\beta_{\gamma}(t_n)$.
Furthermore, for a curve $\gamma \in \curves_{\mathrm{ps}}$, the discrete energy via the trapezoidal rule is defined by
\[
    \energy_{\mathrm{tra},N}^g(\gamma)
    :=
    \frac{1}{N}\sum_{n=0}^{N-1}
    g_{\frac{\gamma(t_n)+\gamma(t_{n+1})}{2}}
    \bigl(
        \beta_{\gamma}(t_n),\beta_{\gamma}(t_n)
    \bigr).
\]
Similarly, the discrete energy via the left-endpoint rule is defined by
\[
    \energy_{\mathrm{l},N}^g(\gamma)
    :=
    \frac{1}{N}\sum_{n=0}^{N-1}
    g_{\gamma(t_n)}
    \bigl(
        \beta_{\gamma}(t_n),\beta_{\gamma}(t_n)
    \bigr).
\]

The error between finite differences and derivatives can be evaluated using the following inequality.
\begin{proposition}
    \label{prop:approx_numerical_integral_beta-gamma}
    Let $\gamma \in \curves_{\mathrm{ps}}$ satisfy
    $\| \dot{\gamma} \|_{L^2}, \| \ddot{\gamma} \|_{L^2} < \infty$.
    Then, for any $t \in [t_n, t_{n+1}]$, we have
    \[
        \| \beta_{\gamma}(t_n) - \dot{\gamma}(t) \|_2
        \leq 2 N^{-1/2} \| \ddot{\gamma} \|_{L^2([t_n, t_{n+1}])}.
    \]
\end{proposition}

\begin{proof}
    By the triangle inequality,
    \[
        \| \beta_{\gamma}(t_n ) - \dot{\gamma} (t) \|_2
        \leq \| \beta_{\gamma}(t_n ) - \dot{\gamma} (t_n) \|_2
        +
        \| \dot{\gamma} (t_n) - \dot{\gamma} (t) \|_2
        .
    \]
    The first term, by Taylor's theorem, equals
    \[
        \beta_{\gamma}(t_n ) - \dot{\gamma} (t_n)
        =
        \frac{1}{h} \int_0^h (h-s)\ddot{\gamma}(t_n+s)\,ds.
    \]
    Thus by the Cauchy--Schwarz inequality,
    \begin{align*}
        \| \beta_{\gamma}(t_n ) - \dot{\gamma} (t_n) \|_2
        &\leq
        \frac{1}{h}
        \|h-s\|_{L^2([0,h])}
        \|\ddot{\gamma}(t_n+\bullet)\|_{L^2([0,h])} \\
        &=
        \frac{h^{1/2}}{\sqrt{3}}
        \|\ddot{\gamma}\|_{L^2([t_n,t_{n+1}])}
        \leq
        N^{-1/2}\|\ddot{\gamma}\|_{L^2([t_n,t_{n+1}])}.
    \end{align*}
    The second term, by applying \Cref{prop:morrey_ineq} to $\dot{\gamma}$, satisfies
    \[
        \| \dot{\gamma} (t_n) - \dot{\gamma} (t) \|_2
        \leq
        N^{-1/2}\|\ddot{\gamma}\|_{L^2([t_n,t_{n+1}])}.
    \]
    Combining these yields the proposition.
\end{proof}

\subsection{Error between curve energy and energy approximated by trapezoidal rule}
\label{sec:error_of_energy_and_trapezoidal_rule}

For constants $K_1, K_2 > 0$, we define the following space of curves $\curves_{K_1, K_2} \subset \curves_{\mathrm{ps}}$:
\[
    \curves_{K_1, K_2} := \left\{ \gamma \in \curves_{\mathrm{ps}}
    \, \middle| \, 
        \left\| \dot{\gamma} \right\|_{L^2([0,1])} \leq K_1 , \,  
        \left\| \ddot{\gamma} \right\|_{L^2([0,1])} \leq K_2   
    \right\}.
\]
In this section, we estimate the error between the curve energy and the energy approximated using the trapezoidal rule for the space of curves $\curves_{K_1, K_2}$.

\begin{proposition}
    \label{prop:functional_lipschitz_estimate}
    Let $H\colon \mathbb{R}^D \to \mathbb{R}^{D\times D}$ be a matrix-valued function such that there exists a constant $c_2 > 0$ satisfying
    \[
        u^{\top} H(x_1) u \le c_2 \|u\|_2^2
        \quad \text{for all } x,u \in \mathbb{R}^D.
    \]
    Suppose that $\| H (x) - H ( y) \|_{\mathcal{B}} \leq L_H \| x - y \|_2$. Then, we have
    \begin{align*}
        & \left|
            u_1^{\top} H(x_1) u_1
            - 
            u_2^{\top} H(x_2) u_2 
        \right|
        \\
        & \leq L_H \| x_1 - x_2 \|_2 \| u_1 \|_2^2  
        + c_2 \| u_1 - u_2 \|_2^2 
        + 2 c_2 \| u_1 \|_2 \| u_1 - u_2 \|_2.
    \end{align*}
\end{proposition}

\begin{proof}
    First, we have
    \begin{align}
        & \left|
            u_1^{\top} H(x_1) u_1
            - 
            u_2^{\top} H(x_2) u_2 
        \right|
        \notag \\
        & \leq
        \left| 
            u_1^{\top} (H(x_1) - H(x_2) ) u_1
        \right|
        +
        \left| 
            ( u_2 - u_1 )^{\top} H(x_2 ) ( u_2 - u_1 )
        \right|
        + 2 \left| 
            u_1^{\top} H(x_2 ) ( u_2 - u_1 )
        \right|.
        \notag 
    \end{align}
    For the first term, we have
    \begin{align*}
        & \left| u_1^{\top} (H(x_1) - H(x_2) ) u_1 \right|
        \\
        & \leq \left\| u_1 \right\|_2 \left\| H(x_1) - H(x_2) \right\|_{\mathcal{B}} \left\| u_1 \right\|_2
        \quad \because \text{Cauchy-Schwarz inequality} \\
        & \leq L_H \left\| x_1 - x_2 \right\|_{2} \left\| u_1 \right\|_2^2.
    \end{align*}
    Similarly, for the second and third terms, we obtain:
    \begin{align*}
        \left| 
            ( u_2 - u_1 )^{\top} H (x_2 ) ( u_2 - u_1 )
        \right|
        & \leq 
        c_2 \| u_1 - u_2 \|_2^2,
        \\
        \left| 
            u_1^{\top} H(x_2 ) ( u_2 - u_1 )
        \right|
        & \leq c_2 \| u_1 \|_2 \| u_1 - u_2 \|_2.
    \end{align*}
\end{proof}

\begin{proposition}
    \label{prop:approx_numerical_integral}
    Let $g$ be a Riemannian metric on $\mathbb{R}^D$ such that, when expressed as $g_x(u,u) = u^{\top} H(x) u$, there exists a constant $c_2 > 0$ satisfying
    \[
        g_x(u,u) \le c_2 \|u\|_2^2
        \quad\text{for all } x,u \in \mathbb{R}^D.
    \]
    Suppose that $\| H (x) - H ( y) \|_{\mathcal{B}} \leq L_H \| x - y \|_2$. Let $\gamma \in \curves_{K_1, K_2}$. Then, we have
    \begin{equation*}
        \left|
            \energy^g (\gamma) -
            \energy_{\mathrm{tra},N}^g (\gamma)
        \right|
        \leq
        \frac{L_H (K_1^3 + 2 K_1^2 K_2) }{N}
        +
        \frac{ 4 c_2  K_1 K_2}{N}
        +
        \frac{4 c_2 K_2^2}{N^2}.
    \end{equation*}
\end{proposition}

\begin{proof}
    Starting from the left-hand side, we have
    \begin{align}
        & \left|
            \energy^g (\gamma) -
            \energy_{\mathrm{tra},N}^g (\gamma)
        \right|
        \notag \\
        & \leq
        \sum_{n=0}^{N-1} \int_{t_n}^{t_{n+1}}
        \biggl| \dot{\gamma} (t)^{\top} H\left(\gamma  (t) \right) \dot{\gamma} (t)
        - \beta_{\gamma} (t_n)^{\top} H \left( \frac{\gamma(t_n) + \gamma (t_{n+1})}{2} \right) \beta_{\gamma} (t_n) \biggr| dt
        \notag \\
        & \leq \sum_{n=0}^{N-1}
        \int_{t_n}^{t_{n+1}} \biggl( L_H \left\| \gamma(t) - \frac{\gamma(t_n) + \gamma (t_{n+1})}{2} \right\|_2 \left\| \dot{\gamma} (t) \right\|_2^2
        \notag \\
        & \quad \hspace{60pt}
        + c_2 \| \dot{\gamma} (t) - \beta_{\gamma} (t_n) \|_2^2
        +
         2 c_2 \| \dot{\gamma} (t) \|_2 \| \dot{\gamma} (t) - \beta_{\gamma} (t_n) \|_2 \biggr) dt
        &
        \because \text{\Cref{prop:functional_lipschitz_estimate}}
        \label{prop:approx_numerical_integral_integrated}.
    \end{align}

    We first estimate the first term in \Cref{prop:approx_numerical_integral_integrated}. By \Cref{prop:morrey_ineq_cor}, for $t \in [t_n , t_{n+1}]$, we have
    \begin{align*}
        \left\| \gamma(t) - \frac{\gamma(t_n) + \gamma (t_{n+1})}{2} \right\|_2
        & \leq 
        \frac{1}{2}
        \left\| \gamma(t) - \gamma(t_n) \right\|_2
        +        
        \frac{1}{2}
        \left\| \gamma(t) - \gamma(t_{n+1}) \right\|_2
        \\
        & \leq 
        \left( \| \dot{\gamma} \|_{L^2([t_n,t_{n+1}])}
        + 2 \| \ddot{\gamma} \|_{L^2([t_n,t_{n+1}])} \right) \frac{1}{N}
        .
    \end{align*}
    Applying this to the first term of \Cref{prop:approx_numerical_integral_integrated}, we obtain
    \begin{align}
        & \sum_{n=0}^{N-1} \int_{t_n}^{t_{n+1}}
        L_H \left\| \gamma(t) - \frac{\gamma(t_n) + \gamma (t_{n+1})}{2} \right\|_2 \left\| \dot{\gamma} (t) \right\|_2^2
        dt 
        \notag \\
        & \leq 
        L_H N^{-1} \sum_{n=0}^{N-1} \left\| \dot{\gamma} (t) \right\|_{L^2 ([t_n,t_{n+1}])}^2 
        \left( 
            \| \dot{\gamma} \|_{L^2([t_n,t_{n+1}])}
            + 2 \| \ddot{\gamma} \|_{L^2([t_n,t_{n+1}])} 
        \right)
        \notag \\
        & \leq 
        L_H N^{-1} \sum_{n=0}^{N-1} \left\| \dot{\gamma} (t) \right\|_{L^2 ([t_n,t_{n+1}])}^2 
        \left( 
            \| \dot{\gamma} \|_{L^2([0,1])}
            + 2 \| \ddot{\gamma} \|_{L^2([0,1])} 
        \right)
        \notag \\
        & \leq 
        L_H N^{-1} \left\| \dot{\gamma} (t) \right\|_{L^2 ([0,1])}^2 
        \left( 
            \| \dot{\gamma} \|_{L^2([0,1])}
            + 2 \| \ddot{\gamma} \|_{L^2([0,1])} 
        \right)
        \leq
        L_H N^{-1} (K_1^3 + 2 K_2 K_1^2 )
        .
        \label{eq:approx_numerical_integral_first}
    \end{align}

    The second term of \Cref{prop:approx_numerical_integral_integrated} is
    \begin{align}
        \sum_{n=0}^{N-1} \int_{t_n}^{t_{n+1}}
            c_2 \| \dot{\gamma} (t) - \beta_{\gamma} (t_n) \|_2^2 dt
        & \leq \sum_{n=0}^{N-1} \int_{t_n}^{t_{n+1}}
            4 c_2 \frac{1}{N} \| \ddot{\gamma} \|_{L^2 ([t_n , t_{n+1}])}^2 dt
        &
        \because \text{\Cref{prop:approx_numerical_integral_beta-gamma}}
        \notag \\
        & \leq  \sum_{n=0}^{N-1} \frac{1}{N} 4 c_2 \frac{1}{N} \| \ddot{\gamma} \|_{L^2 ([t_n , t_{n+1}])}^2
        = 4 c_2 K_2^2 N^{-2} .
        \label{eq:approx_numerical_integral_second}
    \end{align}

    The third term of \Cref{prop:approx_numerical_integral_integrated} is
    \begin{align}
        & \sum_{n=0}^{N-1} \int_{t_n}^{t_{n+1}}
        2 c_2 \| \dot{\gamma} (t) \|_2 \| \dot{\gamma} (t) - \beta_{\gamma} (t_n) \|_2 dt
        \notag \\
        & \leq \sum_{n=0}^{N-1} \int_{t_n}^{t_{n+1}}
        4 c_2 \| \dot{\gamma} (t) \|_2 N^{-1/2} \| \ddot{\gamma} \|_{L^2 ([t_n , t_{n+1}])} dt
        & \because \text{\Cref{prop:approx_numerical_integral_beta-gamma}}
        \notag \\
        & \leq \sum_{n=0}^{N-1} \| 1 \|_{L^2 ([t_n , t_{n+1}])} 
        4 c_2 \| \dot{\gamma} \|_{L^2 ([t_n , t_{n+1}])} 
        N^{-1/2} \| \ddot{\gamma} \|_{L^2 ([t_n , t_{n+1}])}
        & \because \text{Cauchy-Schwarz inequality} \notag 
        \\
        & = \sum_{n=0}^{N-1} N^{-1/2} 
        4 c_2 \| \dot{\gamma} \|_{L^2 ([t_n , t_{n+1}])} N^{-1/2} \| \ddot{\gamma} \|_{L^2 ([t_n , t_{n+1}])}
        \notag \\ 
        & \leq
        4 c_2  K_1 K_2 / N & \because
        \text{Cauchy-Schwarz inequality}.
        \label{eq:approx_numerical_integral_thrid}
    \end{align}

    From \Cref{eq:approx_numerical_integral_first,eq:approx_numerical_integral_second,eq:approx_numerical_integral_thrid}, we get the proposition.
\end{proof}

Hereafter, we set
\[
    C_{\mathrm{tra}}(N)
    :=
    \frac{L_H (K_1^3 + 2 K_1^2 K_2) + 4 c_2 K_1 K_2}{N}
    +
    \frac{4 c_2 K_2^2}{N^2}.
\]
In particular, if $N \ge 1$, then
\[
    C_{\mathrm{tra}}(N)
    \le
    \frac{C_{\mathrm{tra}}^\sharp}{N},
    \quad
    C_{\mathrm{tra}}^\sharp
    :=
    L_H (K_1^3 + 2 K_1^2 K_2) + 4 c_2 K_1 K_2 + 4 c_2 K_2^2.
\]

\begin{proposition}
    \label{prop:approx_numerical_integral_set}
    Under the same setting as \Cref{prop:approx_numerical_integral}, assume $\curves \subset \curves_{K_1, K_2}$. Then
    \[
        \inf_{\gamma \in \curves } \energy_{\mathrm{tra},N}^g ( \gamma )
        \leq
        \inf_{\gamma \in \curves } \energy^g ( \gamma )
        +
        C_{\mathrm{tra}}(N).
    \]
\end{proposition}

\begin{proof}
    For any $\gamma \in \curves$, \Cref{prop:approx_numerical_integral} gives
    \[
        \energy_{\mathrm{tra},N}^g(\gamma)
        \le
        \energy^g(\gamma)+C_{\mathrm{tra}}(N).
    \]
    Taking $\inf_{\gamma \in \curves}$ on both sides yields the conclusion.
\end{proof}

\subsection{Discrete minimizers for the trapezoidal rule and their linear interpolations}
\label{sec:trapolized_and_interpolation}

In this section, we evaluate the error between the discrete energy of a point sequence under the trapezoidal rule and the energy of the curve obtained by linearly interpolating that sequence.

\begin{proposition}
    \label{prop:ineq_cube_3/2}
    For any $x = (x_1 , \ldots , x_N) \in \mathbb{R}^N$, we have
    \[
        \sum_{n=1}^{N} x_n^3 \leq \left( \sum_{n=1}^N x_n^2 \right)^{3/2}.
    \]
\end{proposition}

\begin{proof}
    If the vector $x$ is the zero vector, the proposition holds clearly. 
    Assuming the proposition is correct for the case $\| x \|_2 = 1$, then for any $x \in \mathbb{R}^N \setminus \{ 0 \}$, we have
    \[
        \sum_{n=1}^{N} \left(\frac{x_n}{\| x \|_2}\right)^3 \leq 1,
    \]
    and the proposition follows. Therefore, it suffices to prove the proposition for the case $\| x \|_2 = 1$.

    In the case $\| x \|_2 = 1$, since $|x_n| \leq 1$, we have
    \[
        \sum_{n=1}^{N} x_n^3
        \leq 
        \sum_{n=1}^{N} x_n^2 = 1.
    \]
    Thus, the proposition holds.
\end{proof}

\begin{proposition}
    \label{prop:approx_numerical_integral_pts}
    For any $\gamma^{\mathrm{p}} \in \curves_{\mathrm{p}}^{N}$, let 
    $
        K_3^{2} := \sum_{n=0}^{N-1} \| \beta^{\mathrm{p}} (t_n) \|_2^2 /N.
    $
    Then, we have
    \[
        \left|
            \energy^g (\gamma^{\mathrm{pl}})
            - 
            \energy_{\mathrm{tra},N}^g (\gamma^{\mathrm{p}}) 
        \right|
        \leq
        \frac{L_H K_3^3}{4N^{1/2}}.
    \]
\end{proposition}

\begin{proof}
    For any $0 < s < h$, we have
    \begin{align}
        \gamma^{\mathrm{pl}} (t_n + h/2 ) & = \frac{\gamma^{\mathrm{p}} (t_n) + \gamma^{\mathrm{p}} (t_{n+1} ) }{2},
        \label{eq:pl_midpoint}
        \\
        \dot{\gamma}^{\mathrm{pl}} (t_n + s ) & = \frac{\gamma^{\mathrm{p}} (t_{n+1} ) - \gamma^{\mathrm{p}} (t_n) }{h} = \beta^{\mathrm{p}} (t_n). 
        \label{eq:pl_derivation}   
    \end{align}
    From \Cref{eq:pl_midpoint,eq:pl_derivation}, it follows that 
    \begin{equation}
        \| \dot{\gamma}^{\mathrm{pl}} \|_{L^2 ([0,1])}^2 
        = \frac{1}{N}\sum_{n=0}^{N-1} \| \beta^{\mathrm{p}} (t_n)  \|_2^2
        = K_3^2 .
    \end{equation}

    On the other hand, we have
    \begin{align}
        & 
        \left\| \gamma^{\mathrm{pl}}(t_n + s ) - \frac{\gamma^{\mathrm{p}} (t_n) + \gamma^{\mathrm{p}} (t_{n+1})}{2} \right\|_2
        \notag \\
        & =
        \left\| (1 - Ns ) \gamma^{\mathrm{p}} (t_n) + Ns \gamma^{\mathrm{p}} (t_{n+1}) - \frac{\gamma^{\mathrm{p}} (t_n) + \gamma^{\mathrm{p}} (t_{n+1})}{2} \right\|_2
        \notag \\
        & =
        \left\| \left( \frac{1}{2} - Ns \right) \gamma^{\mathrm{p}} (t_n) + \left( Ns - \frac{1}{2} \right) \gamma^{\mathrm{p}} (t_{n+1}) \right\|_2
        \notag \\
        & =
        \left| Ns - \frac{1}{2} \right|
        \left\| - \gamma^{\mathrm{p}} (t_n) +  \gamma^{\mathrm{p}} (t_{n+1}) \right\|_2
        \notag \\
        & =
        \left| s - \frac{h}{2} \right|
        \left\| \beta^{\mathrm{p}} (t_n) \right\|_2
        \label{eq:pl_traporized}
    \end{align}
    
    From the left-hand side of the proposition, 
    \begin{align}
        & \left| 
            \energy^g (\gamma^{\mathrm{pl}})
            - 
            \energy_{\mathrm{tra},N}^g (\gamma^{\mathrm{p}}) 
        \right|
        \notag \\
        & \leq 
        \sum_{n=0}^{N-1} \int_{t_n}^{t_{n+1}} 
        \biggl| \dot{\gamma}^{\mathrm{pl}} (t)^{\top} H \left(\gamma^{\mathrm{pl}}  (t) \right) \dot{\gamma}^{\mathrm{pl}} (t) 
        \notag \\
        & \hspace{120pt}
        - \beta^{\mathrm{p}} (t_n)^{\top} H \left( \frac{\gamma^{\mathrm{p}}(t_n) + \gamma^{\mathrm{p}} (t_{n+1})}{2} \right) \beta^{\mathrm{p}} (t_n) \biggr| dt
        \notag \\
        & \leq 
        \sum_{n=0}^{N-1} \int_{t_n}^{t_{n+1}} 
        \biggl| \beta^{\mathrm{p}} (t_n)^{\top} H\left(\gamma^{\mathrm{pl}}  (t) \right) \beta^{\mathrm{p}} (t_n) 
        \notag \\
        & \hspace{120pt}
        - \beta^{\mathrm{p}} (t_n)^{\top} H \left( \frac{\gamma^{\mathrm{p}} (t_n) + \gamma^{\mathrm{p}} (t_{n+1})}{2} \right) \beta^{\mathrm{p}} (t_n) \biggr| dt
        & \because \text{\Cref{eq:pl_derivation}}  
        \notag \\
        & \leq \sum_{n=0}^{N-1} 
        \int_{t_n}^{t_{n+1}} 
            L_H 
            \left\| \gamma^{\mathrm{pl}}(t) - \frac{\gamma^{\mathrm{p}} (t_n) + \gamma^{\mathrm{p}} (t_{n+1})}{2} \right\|_2     
            \left\| \beta^{\mathrm{p}} (t_n) \right\|_2^2 dt
        & \because \text{\Cref{prop:functional_lipschitz_estimate}}
        \notag \\
        & = \sum_{n=0}^{N-1} 
        \int_{0}^{h} 
            L_H 
            \left| s - \frac{h}{2} \right|
            \left\| \beta^{\mathrm{p}} (t_n) \right\|_2
            \left\| \beta^{\mathrm{p}} (t_n) \right\|_2^2 ds
        &
        \because \text{\Cref{eq:pl_traporized}}
        \notag \\
        & = \sum_{n=0}^{N-1} 
            \frac{L_H}{4N^2}
            \left\| \beta^{\mathrm{p}} (t_n) \right\|_2^3
        \leq  
            \frac{L_H}{4N^2}
            \left(
                \sum_{n=0}^{N-1}
                \left\| \beta^{\mathrm{p}} (t_n) \right\|_2^2
            \right)^{3/2}
        & \because \text{\Cref{prop:ineq_cube_3/2}}
        \notag \\
        & =
            \frac{L_H}{4} 
            K_3^3
            N^{-1/2}.
        \notag
    \end{align}
\end{proof}

\subsection{Estimation of Christoffel symbols}
\label{sec:estimatation_Christoffel_symbol}

In this section, we evaluate the first derivatives of $g$ and the Christoffel symbols.
\begin{proposition}
    \label{prop:estimate_derivation_g}
    Let $g \in C^1$ be the Riemannian metric in Euclidean space. Suppose that in coordinates $x=(x_1, \ldots , x_{D})$, we can write $g_x (u,u) = u^{\top} H(x) u$ and that $\| H (x) - H ( y) \|_{\mathcal{B}} \leq L_H \| x - y \|_2$. Then, we have
    \[
         \left| \frac{\partial g_{i j}}{\partial x_\ell} \right|
         \leq L_H .
    \]
\end{proposition}
\begin{proof}
    Let $\delta_{ij}$ denote the Kronecker delta, and let $\mathbf{e}_{\ell} := (\delta_{\ell i})_i$, identifying $\mathbf{e}_{\ell}$ with $\partial / \partial x_{\ell}$. Differentiating the function $g_{ij}$ with respect to $\partial / \partial x_{\ell}$, we obtain
    \begin{align}
        \left| \frac{\partial g_{i j}}{\partial x_\ell} \right|
        & =
        \left|
        \lim_{h \to 0 } \frac{1}{h} \left( 
            g_{x + h \mathbf{e}_{\ell}} 
            \left( 
                \frac{\partial}{\partial x_i} , 
                \frac{\partial}{\partial x_j} 
            \right)
            - g_x \left( 
                \frac{\partial}{\partial x_i} , 
                \frac{\partial}{\partial x_j} 
            \right)    
        \right)
        \right|
        \notag \\
        & =
        \left|
        \lim_{h \to 0 } \frac{1}{h} \mathbf{e}_i^{\top} \left( 
             H (x + h \mathbf{e}_{\ell})
             - H (x) 
        \right)
        \mathbf{e}_j
        \right|
        \notag \\
        & =
        \lim_{h \to 0 }
        \left|
         \frac{1}{h} \mathbf{e}_i^{\top} \left( 
             H (x + h \mathbf{e}_{\ell})
             - H (x)
        \right)
        \mathbf{e}_j
        \right|
        \notag \\
        & \leq 
        \lim_{h \to 0 }
        \frac{1}{h}
        \left\| 
             H (x + h \mathbf{e}_{\ell})
             - H (x)
        \right\|_{\mathcal{B}}
        \leq \frac{1}{h} L_H h \| \mathbf{e}_{\ell} \|_2 = L_H.
        \notag 
    \end{align}
\end{proof}

\begin{proposition}
    \label{prop:Christoffel_symbol_estimate}
    Let $g \in C^1$ be a Riemannian metric on $\mathbb{R}^D$. Suppose that we can write $g_x (u,u) = u^{\top} H(x) u$ and assume
    \[
        c_1 \| u \|_2^2 \leq g_x (u,u),
        \quad
        \| H (x) - H ( y) \|_{\mathcal{B}} \leq L_H \| x - y \|_2.
    \]
    Then
    \[
        \left| \Gamma_{ij}^k (x) \right|
        \leq \frac{3 L_H D^{1/2}}{2 c_1}.
    \]
\end{proposition}

\begin{proof}
    By \Cref{prop:estimate_derivation_g}, for any $i,j,\ell$,
    \[
        \left| \frac{\partial g_{ij}}{\partial x_\ell}(x) \right|
        \leq L_H.
    \]
    Hence, by \Cref{eq:Christoffel_symbol},
    \begin{align*}
        |\Gamma_{ij}^k(x)|
        &\leq
        \frac12
        \sum_{\ell=1}^D
        |g^{k\ell}(x)|
        \left(
            \left|\frac{\partial g_{\ell j}}{\partial x_i}(x)\right|
            +
            \left|\frac{\partial g_{i\ell}}{\partial x_j}(x)\right|
            +
            \left|\frac{\partial g_{ij}}{\partial x_\ell}(x)\right|
        \right) \\
        &\leq
        \frac{3 L_H}{2}
        \sum_{\ell=1}^D |g^{k\ell}(x)|.
    \end{align*}

    Set $G(x) := (g_{ij}(x))_{i,j} = H(x)$. From the assumption $c_1 \| u \|_2^2 \le u^\top H(x) u$, we obtain
    \[
        \|G(x)^{-1}\|_{\mathcal B} \le \frac{1}{c_1}.
    \]
    Moreover,
    \begin{align*}
        \sum_{\ell=1}^D |g^{k\ell}(x)|
        &\le
        D^{1/2}
        \left(
            \sum_{\ell=1}^D |g^{k\ell}(x)|^2
        \right)^{1/2} \\
        &=
        D^{1/2}\| e_k^\top G(x)^{-1} \|_2
        \le
        D^{1/2}\|G(x)^{-1}\|_{\mathcal B}
        \le
        \frac{D^{1/2}}{c_1}.
    \end{align*}
    Therefore
    \[
        |\Gamma_{ij}^k(x)|
        \le
        \frac{3 L_H}{2} \cdot \frac{D^{1/2}}{c_1}
        =
        \frac{3 L_H D^{1/2}}{2 c_1}.
    \]
\end{proof}

\subsection{A Priori estimates for minimal geodesics}
\label{sec:apriori_estimate_for_minimal_geodesic}

We derive $L^2$ estimates for the first and second derivatives of minimal geodesics. Define the constants $K_{1,g}$ and $K_{2,g}$ as follows:
\begin{align*}
    K_{1,g} := c_1^{-1/2} d^g ( X^{(0)} , X^{(1)}),
    \quad
    K_{2,g} := \frac{3 L_H D^{2}}{2 c_1} K_{1,g}^2.
\end{align*}

\begin{proposition}
    \label{prop:apriori_estuimate_geodesic}
    Let $g$ be the Riemannian metric on $\mathbb{R}^D$. Suppose that we can express $g_x (u,u) = u^{\top} H(x) u$ and assume the following conditions: $c_1 \| u \|_2^2 \leq g_x (u,u) \leq c_2 \| u \|_2^2$ and $\| H (x) - H ( y) \|_{\mathcal{B}} \leq L_H \| x - y \|_2$. Then for any $\gamma \in \argmin_{\gamma \in \curves_{\mathrm{ps}} (X^{(0)}, X^{(1)})} \energy^g ( \gamma)$, we have
    \[
        \| \dot{\gamma} \|_{L^2} \leq  K_{1,g},
        \quad
        \| \ddot{\gamma} \|_{L^2} 
        \leq 
        K_{2,g}.
    \]
\end{proposition}

\begin{proof}
    Let $\gamma \in \argmin_{\gamma \in \curves_{\mathrm{ps}} (X^{(0)}, X^{(1)})} \energy^g ( \gamma)$.
    By \Cref{prop:bounded_g_completeness,prop:geodesic_len_and_energy_2}, $\gamma$ is a minimal geodesic. From \Cref{prop:metric_geodesic_constant,prop:geodesic_len_and_energy_2}, for any $t \in [0,1]$,
    \[
        g_{\gamma(t)}(\dot{\gamma}(t),\dot{\gamma}(t))
        =
        d^g(X^{(0)},X^{(1)})^2.
    \]
    Hence
    \[
        c_1 \|\dot{\gamma}(t)\|_2^2
        \le
        d^g(X^{(0)},X^{(1)})^2,
    \]
    so
    \[
        \|\dot{\gamma}(t)\|_2
        \le
        c_1^{-1/2} d^g(X^{(0)},X^{(1)})
        =
        K_{1,g}
        \quad (t\in[0,1]).
    \]
    In particular,
    \[
        \|\dot{\gamma}\|_{L^2} \le K_{1,g}.
    \]

    Next, $\gamma$ satisfies the geodesic equation
    \[
        \ddot{\gamma}_k
        =
        -\sum_{i,j=1}^D \Gamma_{ij}^k(\gamma)\dot{\gamma}_i\dot{\gamma}_j,
    \]
    so
    \begin{align*}
        |\ddot{\gamma}_k(t)|
        &\le
        \sum_{i,j=1}^D
        |\Gamma_{ij}^k(\gamma(t))|
        |\dot{\gamma}_i(t)|
        |\dot{\gamma}_j(t)| \\
        &\le
        \sup_{x,i,j,k} |\Gamma_{ij}^k(x)|
        \left( \sum_{i=1}^D |\dot{\gamma}_i(t)| \right)^2 \\
        &\le
        \sup_{x,i,j,k} |\Gamma_{ij}^k(x)|\,
        D \|\dot{\gamma}(t)\|_2^2
        \quad \text{(Cauchy--Schwarz inequality)} \\
        &\le
        \frac{3L_H D^{1/2}}{2c_1}\,D\,K_{1,g}^2
        =
        \frac{3L_H D^{3/2}}{2c_1} K_{1,g}^2
        \quad \text{by \Cref{prop:Christoffel_symbol_estimate}}.
    \end{align*}
    Therefore
    \[
        \|\ddot{\gamma}(t)\|_2
        \le
        D^{1/2}\max_{1\le k\le D} |\ddot{\gamma}_k(t)|
        \le
        \frac{3L_H D^2}{2c_1}K_{1,g}^2
        =
        K_{2,g}.
    \]
    Consequently
    \[
        \|\ddot{\gamma}\|_{L^2} \le K_{2,g}.
    \]
\end{proof}

\subsection{Proof of Theorem \ref{theo:approx_numerical_integral_trapezoidal}}
\label{sec:proof_of_theo:approx_numerical_integral_trapezoidal}

\begin{proof}[Proof of \Cref{theo:approx_numerical_integral_trapezoidal}]
    From \Cref{prop:geodesic_len_and_energy_2},
    \[
        \min_{\gamma \in \curves_{\mathrm{geo}} (X^{(0)}, X^{(1)}) }
        \energy^g ( \gamma )
        =
        \min_{\gamma \in \curves_{\mathrm{ps}} (X^{(0)}, X^{(1)}) }
        \energy^g ( \gamma ).
    \]
        Since
    \[
        \curves_{\mathrm{geo}} (X^{(0)}, X^{(1)})
        \subset
        \curves_{\mathrm{ps}} (X^{(0)}, X^{(1)})
    \]
    and
    \[
        \min_{\gamma \in \curves_{\mathrm{geo}} (X^{(0)}, X^{(1)}) }
        \energy^g ( \gamma )
        =
        \min_{\gamma \in \curves_{\mathrm{ps}} (X^{(0)}, X^{(1)}) }
        \energy^g ( \gamma ),
    \]
    we have
    \[
        \argmin_{\gamma \in \curves_{\mathrm{geo}} (X^{(0)}, X^{(1)}) }
        \energy^g ( \gamma )
        \subset
        \argmin_{\gamma \in \curves_{\mathrm{ps}} (X^{(0)}, X^{(1)}) }
        \energy^g ( \gamma ).
    \]
    Hence, by \Cref{prop:apriori_estuimate_geodesic},
    \[
        \argmin_{\gamma \in \curves_{\mathrm{geo}} (X^{(0)}, X^{(1)}) }
        \energy^g ( \gamma )
        \subset
        \curves_{K_{1,g},K_{2,g}}.
    \]
        Set
    \[
        \mathcal{A}
        :=
        \argmin_{\gamma \in \curves_{\mathrm{geo}} (X^{(0)}, X^{(1)}) }
        \energy^g ( \gamma ).
    \]
    By \Cref{prop:geodesic_len_and_energy_2}, the minimum energy over
    $\curves_{\mathrm{ps}} (X^{(0)}, X^{(1)})$ is attained by a geodesic.
    Hence $\mathcal{A} \neq \emptyset$.
    Applying \Cref{prop:approx_numerical_integral_set} with $\curves = \mathcal{A}$, we obtain
    \begin{align}
        \min_{\gamma^{\mathrm{p}} \in \curves_{\mathrm{p}}^{N} (X^{(0)}, X^{(1)}) }
        \energy_{\mathrm{tra},N}^g ( \gamma^{\mathrm{p}} )
        &\le
        \inf_{\gamma \in \mathcal{A}}
        \energy_{\mathrm{tra},N}^g ( \gamma ) \notag\\
        &\le
        \inf_{\gamma \in \mathcal{A}}
        \energy^g ( \gamma )
        +
        C_{\mathrm{tra}}(N) \notag\\
        &=
        \min_{\gamma \in \curves_{\mathrm{geo}} (X^{(0)}, X^{(1)}) }
        \energy^g ( \gamma )
        +
        C_{\mathrm{tra}}(N) \notag\\
        &=
        \min_{\gamma \in \curves_{\mathrm{ps}} (X^{(0)}, X^{(1)}) }
        \energy^g ( \gamma )
        +
        C_{\mathrm{tra}}(N).
        \label{eq:approx_numerical_integral_lower_bound_revised}
    \end{align}
    Hence
    \[
        \energy_{\mathrm{tra},N}^g(\gamma_{\mathrm{tra},N}^{*\mathrm{p}})
        \le
        \min_{\gamma \in \curves_{\mathrm{ps}} (X^{(0)}, X^{(1)}) }
        \energy^g(\gamma)
        +
        \frac{C_{\mathrm{tra}}^\sharp}{N}.
    \]

    Furthermore,
    \[
        \frac1N\sum_{n=0}^{N-1}\|\beta^{\mathrm{p}}(t_n)\|_2^2
        \le
        \frac{\energy_{\mathrm{tra},N}^g(\gamma^{\mathrm{p}})}{c_1},
    \]
    so for $\gamma^{\mathrm{p}} = \gamma_{\mathrm{tra},N}^{*\mathrm{p}}$,
    \[
        \frac1N\sum_{n=0}^{N-1}\|\beta^{\mathrm{p}}(t_n)\|_2^2
        \le
        \frac{ d^g(X^{(0)},X^{(1)})^2 + C_{\mathrm{tra}}^\sharp}{c_1}
        =: \overline K_{3,g}^{\,2}.
    \]
    By \Cref{prop:approx_numerical_integral_pts},
    \[
        \energy^g(\gamma_{\mathrm{tra},N}^{*\mathrm{pl}})
        \le
        \energy_{\mathrm{tra},N}^g(\gamma_{\mathrm{tra},N}^{*\mathrm{p}})
        +
        \frac{L_H \overline K_{3,g}^{\,3}}{4N^{1/2}}.
    \]
    Therefore
    \begin{align*}
        \min_{\gamma \in \curves_{\mathrm{ps}} (X^{(0)}, X^{(1)}) }
        \energy^g(\gamma)
        &\le
        \energy^g(\gamma_{\mathrm{tra},N}^{*\mathrm{pl}}) \\
        &\le
        \energy_{\mathrm{tra},N}^g(\gamma_{\mathrm{tra},N}^{*\mathrm{p}})
        +
        \frac{L_H \overline K_{3,g}^{\,3}}{4N^{1/2}} \\
        &\le
        \min_{\gamma \in \curves_{\mathrm{ps}} (X^{(0)}, X^{(1)}) }
        \energy^g(\gamma)
        +
        \frac{C_{\mathrm{tra}}^\sharp}{N}
        +
        \frac{L_H \overline K_{3,g}^{\,3}}{4N^{1/2}}.
    \end{align*}
    Finally, since $N^{-1} \le N^{-1/2}$ for $N \ge 1$, setting
    \[
        C
        :=
        C_{\mathrm{tra}}^\sharp
        +
        \frac{L_H \overline K_{3,g}^{\,3}}{4}
    \]
    yields
    \[
        \energy_{\mathrm{tra},N}^g(\gamma_{\mathrm{tra},N}^{*\mathrm{p}})
        \le
        \min_{\gamma \in \curves_{\mathrm{ps}} (X^{(0)}, X^{(1)}) }
        \energy^g(\gamma)
        +
        \frac{C}{N^{1/2}},
    \]
    and
    \[
        \energy^g(\gamma_{\mathrm{tra},N}^{*\mathrm{pl}})
        \le
        \energy_{\mathrm{tra},N}^g(\gamma_{\mathrm{tra},N}^{*\mathrm{p}})
        +
        \frac{C}{N^{1/2}}.
    \]
    This gives the theorem.
\end{proof}

\section{Approximation of energy minimization problem via left-endpoint rule}
\label{sec:proof_of_approx_trapezoidal_left}

In this section, we prove \Cref{theo:approx_numerical_integral_left}. In \Cref{sec:leftendpoint_and_trapezoidal_rule}, we compare the energy approximated by the left-endpoint rule with that approximated by the trapezoidal rule. In \Cref{sec:leftendpopint_and_interpolation}, we evaluate the error between the energy approximated for points obtained using the left-endpoint rule and the energy of a curve formed by linear interpolation of these points. Finally, in \Cref{sec:proof_of_theo:approx_numerical_integral_left}, we prove \Cref{theo:approx_numerical_integral_left}.

\subsection{Left-endpoint rule and trapezoidal rule}
\label{sec:leftendpoint_and_trapezoidal_rule}

In this section, we compare the energy approximated by the left-endpoint rule with that approximated by the trapezoidal rule.
\begin{proposition}
    \label{prop:left_trapezoidal_fomula}
    Let $g \in C^1$ be the Riemannian metric in $\mathbb{R}^D$. Suppose that we can express $g_x (u,u) = u^{\top} H(x) u$ and that $\| H(x) - H(y) \|_{\mathcal{B}} \leq L_H \| x - y \|_2$. For any $\gamma^{\mathrm{p}} \in \curves_{\mathrm{p}}^{N}$, let 
    \[
        K_3^{2} := \frac{1}{N} \sum_{n=0}^{N-1} \| \beta^{\mathrm{p}} (t_n) \|_2^2.
    \]
    Then,
    \[
        \left|
            \energy_{\mathrm{l},N}^g (\gamma^{\mathrm{p}})
            -
            \energy_{\mathrm{tra},N}^g (\gamma^{\mathrm{p}})
        \right|
        \leq
        \frac{L_H K_3^3}{2N^{1/2}}.
    \]
\end{proposition}

\begin{proof}
    From the left-hand side of the proposition, 
    \begin{align}
        & \left| 
            \energy_{\mathrm{l},N}^g (\gamma^{\mathrm{p}})
            - 
            \energy_{\mathrm{tra},N}^g (\gamma^{\mathrm{p}}) 
        \right|
        \notag \\
        & \leq 
        \sum_{n=0}^{N-1} 
        \biggl| \beta^{\mathrm{p}} (t_n)^{\top} H\left(\gamma^{\mathrm{p}}  (t_n) \right) \beta^{\mathrm{p}} (t_n) 
        - \beta^{\mathrm{p}} (t_n)^{\top} H \left( \frac{\gamma^{\mathrm{p}} (t_n) + \gamma^{\mathrm{p}} (t_{n+1})}{2} \right) \beta^{\mathrm{p}} (t_n) \biggr|
        & \because \text{\Cref{eq:pl_derivation}}  
        \notag \\
        & \leq 
        \frac{1}{N} \sum_{n=0}^{N-1} 
            L_H 
            \left\| \gamma^{\mathrm{p}}(t_n) - \frac{\gamma^{\mathrm{p}} (t_n) + \gamma^{\mathrm{p}} (t_{n+1})}{2} \right\|_2     
            \left\| \beta^{\mathrm{p}} (t_n) \right\|_2^2
        & \because \text{\Cref{prop:functional_lipschitz_estimate}}
        \notag \\
        & = 
        \frac{L_H}{2 N^2} \sum_{n=0}^{N-1} 
            \left\| \beta^{\mathrm{p}} (t_n) \right\|_2^3
        \leq  
            \frac{L_H}{2N^2}
            \left(
                \sum_{n=0}^{N-1}
                \left\| \beta^{\mathrm{p}} (t_n) \right\|_2^2
            \right)^{3/2}
        & \because \text{\Cref{prop:ineq_cube_3/2}}
        \notag \\
        & =
            \frac{L_H}{2} 
            K_3^3
            N^{-1/2}
            .
            \notag
    \end{align}
\end{proof}

\subsection{Discrete minimizers for the left-endpoint rule and their linear interpolations}
\label{sec:leftendpopint_and_interpolation}

In this section, we estimate the error between the left-endpoint discrete energy of a point sequence and the energy of the curve obtained by linearly interpolating that sequence.

\begin{proposition}
    \label{prop:approx_numerical_integral_pts_left}
    Let $g \in C^1$ be the Riemannian metric on $\mathbb{R}^D$. Suppose that we can express $g_x (u,u) = u^{\top} H(x) u$ and that $\| H(x) - H(y) \|_{\mathcal{B}} \leq L_H \| x - y \|_2$. For any $\gamma^{\mathrm{p}} \in \curves_{\mathrm{p}}^{N}$, let
    \[
        K_3^{2} := \frac{1}{N} \sum_{n=0}^{N-1} \| \beta^{\mathrm{p}} (t_n) \|_2^2.
    \]
    Then,
    \[
        \left|
            \energy^g (\gamma^{\mathrm{pl}})
            - 
            \energy_{\mathrm{l},N}^g (\gamma^{\mathrm{p}})
        \right|
        \leq
        \frac{L_H K_3^3}{2N^{1/2}}.
    \]
\end{proposition}

\begin{proof}
    For $0 \leq s \leq h$,
    \begin{align}
        & 
        \left\| \gamma^{\mathrm{pl}}(t_n + s ) - \gamma^{\mathrm{p}} (t_n) \right\|_2
        \notag \\
        & =
        \left\| (1 - Ns ) \gamma^{\mathrm{p}} (t_n) + Ns \gamma^{\mathrm{p}} (t_{n+1}) - \gamma^{\mathrm{p}} (t_n) \right\|_2
        \notag \\
        & =
        \left\|  - Ns \gamma^{\mathrm{p}} (t_n) + Ns \gamma^{\mathrm{p}} (t_{n+1}) \right\|_2
        \notag \\
        & =
         Ns 
        \left\| - \gamma^{\mathrm{p}} (t_n) +  \gamma^{\mathrm{p}} (t_{n+1}) \right\|_2
        \notag \\
        & =
        s
        \left\| \beta^{\mathrm{p}} (t_n) \right\|_2
        .
        \label{eq:pl_left}
    \end{align}  
    The left-hand side of the proposition is
    \begin{align}
        & \left| 
            \energy^g (\gamma^{\mathrm{pl}})
            - 
            \energy_{\mathrm{l},N}^g (\gamma^{\mathrm{p}}) 
        \right|
        \notag \\
        & \leq 
        \sum_{n=0}^{N-1} \int_{t_n}^{t_{n+1}} 
        \biggl| \dot{\gamma}^{\mathrm{pl}} (t)^{\top} H \left(\gamma^{\mathrm{pl}}  (t) \right) \dot{\gamma}^{\mathrm{pl}} (t) 
        - \beta^{\mathrm{p}} (t_n)^{\top} H \left( \gamma^{\mathrm{p}}(t_n)  \right) \beta^{\mathrm{p}} (t_n) \biggr| dt
        \notag \\
        & \leq 
        \sum_{n=0}^{N-1} \int_{t_n}^{t_{n+1}} 
        \biggl| \beta^{\mathrm{p}} (t_n)^{\top} H\left(\gamma^{\mathrm{pl}}  (t) \right) \beta^{\mathrm{p}} (t_n)
        - \beta^{\mathrm{p}} (t_n)^{\top} H \left( \gamma^{\mathrm{p}}(t_n)  \right) \beta^{\mathrm{p}} (t_n) \biggr| dt
        & \because \text{\Cref{eq:pl_derivation}}  
        \notag \\
        & \leq \sum_{n=0}^{N-1} 
        \int_{t_n}^{t_{n+1}} 
            L_H 
            \left\| \gamma^{\mathrm{pl}}(t ) - \gamma^{\mathrm{p}} (t_n) \right\|_2     
            \left\| \beta^{\mathrm{p}} (t_n) \right\|_2^2 dt
        & \because \text{\Cref{prop:functional_lipschitz_estimate}}
        \notag \\
        & = \sum_{n=0}^{N-1} 
        \int_{0}^{h} 
            L_H 
            s
            \left\| \beta^{\mathrm{p}} (t_n) \right\|_2
            \left\| \beta^{\mathrm{p}} (t_n) \right\|_2^2 ds
        &
        \because \text{\Cref{eq:pl_left}}
        \notag \\
        & = \sum_{n=0}^{N-1} 
            \frac{L_H}{2N^2}
            \left\| \beta^{\mathrm{p}} (t_n) \right\|_2^3
        \leq  
            \frac{L_H}{2N^2}
            \left(
                \sum_{n=0}^{N-1}
                \left\| \beta^{\mathrm{p}} (t_n) \right\|_2^2
            \right)^{3/2}
        & \because \text{\Cref{prop:ineq_cube_3/2}}
        \notag \\
        & =
            \frac{L_H}{2} 
            K_3^3
            N^{-1/2}.
        \notag
    \end{align}
\end{proof}

\subsection{Proof of Theorem \ref{theo:approx_numerical_integral_left}}
\label{sec:proof_of_theo:approx_numerical_integral_left}

\begin{proof}[Proof of \Cref{theo:approx_numerical_integral_left}]
    Use
    \[
        \overline K_{3,g}^{\,2}
        :=
        \frac{ d^g(X^{(0)},X^{(1)})^2 + C_{\mathrm{tra}}^\sharp}{c_1}
    \]
    from the proof of \Cref{theo:approx_numerical_integral_trapezoidal}. By \Cref{prop:left_trapezoidal_fomula}, for any $\gamma^{\mathrm{p}} \in \curves_{\mathrm{p}}^N(X^{(0)},X^{(1)})$ satisfying
    \[
        \frac1N\sum_{n=0}^{N-1}\|\beta^{\mathrm{p}}(t_n)\|_2^2
        \le
        \overline K_{3,g}^{\,2},
    \]
    we have
    \[
        \left|
            \energy_{\mathrm{l},N}^g(\gamma^{\mathrm{p}})
            -
            \energy_{\mathrm{tra},N}^g(\gamma^{\mathrm{p}})
        \right|
        \le
        \frac{L_H \overline K_{3,g}^{\,3}}{2N^{1/2}}.
    \]
    In particular, for $\gamma_{\mathrm{tra},N}^{*\mathrm{p}}$,
    \[
        \energy_{\mathrm{l},N}^g(\gamma_{\mathrm{tra},N}^{*\mathrm{p}})
        \le
        \energy_{\mathrm{tra},N}^g(\gamma_{\mathrm{tra},N}^{*\mathrm{p}})
        +
        \frac{L_H \overline K_{3,g}^{\,3}}{2N^{1/2}}.
    \]
    Hence
    \begin{align}
        \min_{\gamma^{\mathrm{p}} \in \curves_{\mathrm{p}}^N(X^{(0)},X^{(1)})}
        \energy_{\mathrm{l},N}^g(\gamma^{\mathrm{p}})
        &\le
        \energy_{\mathrm{l},N}^g(\gamma_{\mathrm{tra},N}^{*\mathrm{p}}) \notag\\
        &\le
        \energy_{\mathrm{tra},N}^g(\gamma_{\mathrm{tra},N}^{*\mathrm{p}})
        +
        \frac{L_H \overline K_{3,g}^{\,3}}{2N^{1/2}} \notag\\
        &\le
        \min_{\gamma \in \curves_{\mathrm{ps}}(X^{(0)},X^{(1)})}
        \energy^g(\gamma)
        +
        \frac{C_{\mathrm{tra}}^\sharp}{N}
        +
        \frac{L_H \overline K_{3,g}^{\,3}}{2N^{1/2}}.
        \label{eq:left_upper_bound_revised}
    \end{align}

    Now, for
    \[
        \gamma_{\mathrm{l},N}^{*\mathrm{p}}
        \in
        \argmin_{\gamma^{\mathrm{p}} \in \curves_{\mathrm{p}}^N(X^{(0)},X^{(1)})}
        \energy_{\mathrm{l},N}^g(\gamma^{\mathrm{p}}),
    \]
    we have
    \[
        \frac1N\sum_{n=0}^{N-1}\|\beta^{\mathrm{p}}(t_n)\|_2^2
        \le
        \frac{\energy_{\mathrm{l},N}^g(\gamma_{\mathrm{l},N}^{*\mathrm{p}})}{c_1}.
    \]
    From \Cref{eq:left_upper_bound_revised} and $N^{-1} \le N^{-1/2}$,
    \[
        \energy_{\mathrm{l},N}^g(\gamma_{\mathrm{l},N}^{*\mathrm{p}})
        \le
        \min_{\gamma \in \curves_{\mathrm{ps}}(X^{(0)},X^{(1)})}
        \energy^g(\gamma)
        +
        \frac{
            C_{\mathrm{tra}}^\sharp + \frac12 L_H \overline K_{3,g}^{\,3}
        }{N^{1/2}}.
    \]
    Therefore
    \[
        \frac1N\sum_{n=0}^{N-1}\|\beta^{\mathrm{p}}(t_n)\|_2^2
        \le
        \frac{
            d^g(X^{(0)},X^{(1)})^2
            +
            C_{\mathrm{tra}}^\sharp
            +
            \frac12 L_H \overline K_{3,g}^{\,3}
        }{c_1}
        =:
        \overline K_{4,g}^{\,2}.
    \]
    Thus, by \Cref{prop:approx_numerical_integral_pts_left},
    \[
        \energy^g(\gamma_{\mathrm{l},N}^{*\mathrm{pl}})
        \le
        \energy_{\mathrm{l},N}^g(\gamma_{\mathrm{l},N}^{*\mathrm{p}})
        +
        \frac{L_H \overline K_{4,g}^{\,3}}{2N^{1/2}}.
    \]
    Combining the above, there exists a constant $C$ independent of $N$, namely
    \[
        C
        :=
        C_{\mathrm{tra}}^\sharp
        +
        \frac12 L_H \overline K_{3,g}^{\,3}
        +
        \frac12 L_H \overline K_{4,g}^{\,3},
    \]
    such that
    \[
        \energy_{\mathrm{l},N}^g(\gamma_{\mathrm{l},N}^{*\mathrm{p}})
        \le
        \min_{\gamma \in \curves_{\mathrm{ps}}(X^{(0)},X^{(1)})}
        \energy^g(\gamma)
        +
        \frac{C}{N^{1/2}},
    \]
    and
    \[
        \energy^g(\gamma_{\mathrm{l},N}^{*\mathrm{pl}})
        \le
        \energy_{\mathrm{l},N}^g(\gamma_{\mathrm{l},N}^{*\mathrm{p}})
        +
        \frac{C}{N^{1/2}}.
    \]
    This gives the theorem.
\end{proof}

\section{Length minimization and numerical integration}
\label{sec:proof_of_length_approx_trapezoidal_left}

In this section, we prove \Cref{theo:approx_numerical_integral_left_length,theo:approx_numerical_integral_trapezoidal_length}. The length of a curve can be estimated in terms of its energy.

\begin{proposition}
    \label{prop:length_energy_estimation}
    Let $(M , g )$ be a complete Riemannian manifold. Let $X^{(0)}, X^{(1)} \subset M$ be non-empty closed sets, such that at least one is bounded. We define the curve $\gamma^*\in \curves_{\mathrm{ps}}$ by
    \[
        \gamma^* \in \argmin_{\gamma \in \curves_{\mathrm{ps}} (X^{(0)} , X^{(1)} )} \length^g ( \gamma).
    \]
    Then, for any curve $\gamma \in \curves_{\mathrm{ps}}(X^{(0)} , X^{(1)} )$, we have
    \[
        \length^g (\gamma)^2
        -
        \length^g (\gamma^*)^2
        \leq
        \energy^g (\gamma) 
        - 
        \energy^g (\gamma^*) 
    .
    \]
\end{proposition}

\begin{proof}
    In general, from the Cauchy-Schwarz inequality, we have
    \[
        \length^g (\gamma)^2 \leq \energy^g (\gamma). 
    \]
    From \Cref{prop:geodesic_len_and_energy_2}, it follows that
    \[
        \length^g (\gamma)^2
        -
        \length^g (\gamma^*)^2
        \leq 
        \energy^g (\gamma) 
        -
        \energy^g (\gamma^*)
        .
    \]
\end{proof}

\begin{proof}[Proof of \Cref{theo:approx_numerical_integral_trapezoidal_length,theo:approx_numerical_integral_left_length}]
    From \Cref{prop:length_energy_estimation,theo:approx_numerical_integral_trapezoidal}, we obtain \Cref{theo:approx_numerical_integral_trapezoidal_length}. Furthermore, from \Cref{prop:length_energy_estimation,theo:approx_numerical_integral_left}, we obtain \Cref{theo:approx_numerical_integral_left_length}.
\end{proof}

\begin{remark}
    For the left-endpoint discretization, the $O(N^{-1})$ convergence of the discrete minimum value is consistent with the framework of \cite{rumpf2015variational}. The point of the present argument is different: it is aimed at the energies and squared lengths of the linearly interpolated discrete minimizers in the set-to-set setting. For these reconstructed continuous quantities, we obtain the $O(N^{-1/2})$ estimates proved above. Whether these reconstruction rates can be improved under additional discrete regularity or equidistribution information is left open here.
\end{remark}

\section*{Generative AI}
In this paper, we used GPT-5.2, GPT-5.4, and Opus 4.7 for editing the text.

\bibliographystyle{plain}
\bibliography{main_EN} 

\end{document}